\documentclass[11pt,a4paper]{amsart}
\usepackage[T1]{fontenc}
\usepackage{lmodern}
\usepackage[margin=29mm]{geometry}
\usepackage{amssymb,mathtools,mathrsfs}
\usepackage{microtype}
\usepackage{enumitem}
\usepackage{xcolor}
\usepackage[colorlinks=true,linkcolor=blue!45!black,citecolor=blue!45!black,urlcolor=blue!45!black]{hyperref}
\usepackage{bookmark}
\usepackage{centernot}
\usepackage[only,trianglelefteqslant]{stmaryrd}
\let\leq\leqslant

\let\geq\geqslant

\let\trianglelefteq\trianglelefteqslant

\hypersetup{pdftitle={Melnikov residuals and stationary amalgamation in PAC fields},pdfauthor={}}
\setlist[enumerate]{itemsep=2pt,topsep=4pt}
\newtheorem{theorem}{Theorem}[section]
\newtheorem{lemma}[theorem]{Lemma}
\newtheorem{proposition}[theorem]{Proposition}
\newtheorem{fact}[theorem]{Fact}
\newtheorem{corollary}[theorem]{Corollary}
\theoremstyle{definition}
\newtheorem{definition}[theorem]{Definition}
\newtheorem{remark}[theorem]{Remark}
\DeclareMathOperator{\Gal}{Gal}
\DeclareMathOperator{\acl}{acl}
\DeclareMathOperator{\Th}{Th}
\DeclareMathOperator{\Imfin}{Im}
\DeclareMathOperator{\Frac}{Frac}

\newcommand{\D}{\mathcal D}

\newcommand{\K}{\mathscr K}
\newcommand{\U}{\mathbb{U}}
\newcommand{\Fhat}{\widehat F}

\newcommand{\RD}{R_{\D}}
\newcommand{\sep}{\mathrm{sep}}
\def\Ind#1#2{%
  #1\setbox0=\hbox{$#1x$}%
  \kern\wd0
  \hbox to 0pt{\hss$#1\mid$\hss}%
  \lower.9\ht0\hbox to 0pt{\hss$#1\smile$\hss}%
  \kern\wd0
}

\newcommand{\ind}{%
  \mathop{\mathpalette\Ind{}}\displaylimits^{\D}%
}

\newcommand{\indfork}{%
  \mathop{\mathpalette\Ind{}}\displaylimits^{f}%
}

\newcommand{\indKim}{%
  \mathop{\mathpalette\Ind{}}\displaylimits^{K}%
}
\newcommand{\indld}{%
  \mathop{\mathpalette\Ind{}}\displaylimits^{\mathrm{ld}}%
}

\newcommand{\indgen}{%
  \mathop{\mathpalette\Ind{}}\displaylimits%
}

\newcommand{\nind}{%
  \mathop{\centernot{\mathpalette\Ind{}}}\displaylimits^{\D}%
}

\newcommand{\nindfork}{%
  \mathop{\centernot{\mathpalette\Ind{}}}\displaylimits^{f}%
}

\newcommand{\nindKim}{%
  \mathop{\centernot{\mathpalette\Ind{}}}\displaylimits^{K}%
}

\title[Melnikov residuals and PAC amalgamation]{Melnikov residuals and stationary amalgamation in PAC fields}
\author[\"{O}.~Beyarslan]{\"{O}zlem Beyarslan}
\address{Bo\u{g}azi\c{c}i \"{U}niversitesi}
\email{ozlem.beyarslan@boun.edu.tr}
\thanks{The first author acknowledges support from the
T\"{u}bitak 1001 grant no.~124F359}

\author[P.~Kowalski]{Piotr Kowalski}
\address{Instytut Matematyczny\\
Uniwersytet Wroc{\l}awski\\
Wroc{\l}aw\\
Poland}
\email{pkowa@math.uni.wroc.pl}
\urladdr{http://www.math.uni.wroc.pl/~pkowa/}
\thanks{The second author acknowledges support from the
Narodowe Centrum Nauki grant no.~2018/31/B/ST1/00357
and the T\"{u}bitak 1001 grant no.~124F359.}
\date{}

\begin{document}
\begin{abstract}
We construct an invariant symmetric ternary relation with existence
and stationarity over models for classes of characteristic-zero
Frobenius fields associated with Melnikov formations.
For proper formations, the automorphism theories of their
canonical root-predicate expansions have amalgamation and
no model companion. The underlying fields are then unbounded,
and their theories are NSOP$_1$ but not simple.
In this case, on triples of elementary submodels, the constructed
relation implies ordinary nonforking, which implies
Kim-independence, which implies linear disjointness;
all three implications are strict.
We also study p.e.c. fields
with group actions, namely fields existentially closed
in every regular extension carrying an extension of the action.
We prove that the class of p.e.c.\ difference fields
of characteristic zero is not elementary,
and that the theory of all characteristic-zero difference fields
with canonical root predicates has no model companion.
\end{abstract}

\maketitle
\tableofcontents

\section{Introduction}\label{sec:intro}

The theory of algebraically closed fields with an automorphism
has a model companion, ACFA. This raises the question of which
other theories of fields admit model companions after adding
an automorphism. The present paper is a companion to
\cite{BK2026}, which studies this question for bounded PAC fields.
Here we consider characteristic-zero Frobenius fields associated
with Melnikov formations, with particular emphasis on proper
formations. We construct a ternary relation on elementary
submodels and use it to prove amalgamation for the automorphism
theories of the canonical root-predicate expansions.
For proper formations, we prove that these automorphism theories
have no model companion.
The formation of all finite groups gives the familiar algebraically
closed case, whose model companion is $\mathrm{ACFA}_0$ up to
the canonical root-predicate expansion.

The construction uses the residual theory of Fried--Haran
\cite{FHQuasi,FHTheta} for Melnikov formations of finite groups.
The relation combines linear disjointness with a prescribed
relative algebraic closure of the compositum.
We prove that it is invariant and symmetric and has existence
and stationarity over models.

A related question concerns existential closure in regular
extensions of fields with group actions.
A field extension is \emph{regular} if it is separable and
the base field is relatively algebraically closed in it.
Motivated by the study of group actions in
\cite{BKVirtuallyFree,BK}, we ask whether this restricted
closure condition is axiomatizable.
We show that this condition is not elementary for a single
automorphism in characteristic zero, whereas for every finite
group it is elementary in arbitrary characteristic.
These ordinary-language results are proved directly in
Section~\ref{sec:regular-ec}, where we also deduce the nonexistence
of a model companion in the root-predicate language using
Corollary~\ref{cor:noMC}, applied to the canonical root-predicate
expansions of $\omega$-free PAC fields.

For the amalgamation argument, we use the ring language
expanded by predicates expressing that monic polynomials
have a root. For a field $K$, write $K_S$
for its canonical definitional expansion by these predicates,
taken from the Galois language of
Cherlin--van den Dries--Macintyre \cite{CDM};
see Section~\ref{sec:realization}.
In characteristic zero, embeddings between these expansions
are precisely regular field embeddings.
For the Frobenius fields considered here, the complete theories
of these expansions eliminate quantifiers and are therefore
model complete.
For proper formations, these theories also have the independence
property and are therefore unstable.
Our amalgamation theorem then supplies the remaining hypothesis
of Kikyo's theorem \cite[Theorem~3.3]{Kikyo}, which rules out
a model companion for each corresponding automorphism theory.

We isolate a general criterion for the amalgamation
property for automorphisms (PAPA): invariance, existence,
and stationarity of a ternary relation on small elementary
submodels of a complete theory suffice
(Lemma~\ref{lem:papa-criterion}).
The new ternary relation constructed on the PAC fields
considered here provides instances of this criterion in
nonsimple theories. For proper formations, it is strictly
stronger than ordinary nonforking and Kim-independence
(Section~\ref{sec:kim-comparison}).
We also recall the general ultraproduct criterion for
elementary axiomatizability of existentially closed models:
for any first-order theory, the class of all its e.c.\ models
is elementary if and only if it is closed under ultraproducts;
for inductive theories, these conditions are equivalent to
the existence of a model companion
(Fact~\ref{fact:ec-ultraproducts};
see also \cite[Introduction]{GSTEcActions}).
For the inductive root-predicate theory, this connects
nonexistence of a model companion with failure of
ultraproduct closure of its e.c.\ models.

We introduce the notation needed to state the main result.
We call a nonempty class $\D$ of finite groups a
\emph{Melnikov formation} if it is closed under isomorphisms,
normal subgroups, quotients, and extensions.
We call $\D$ \emph{proper} if it is different from the class
of all finite groups.

For a field $E$, let $E^{\mathrm{alg}}$ denote an algebraic
closure of $E$, and let $E^{\sep}$ denote its subfield
consisting of the elements separable over $E$.
Throughout the paper, these closures are chosen compatibly:
for subfields of a common field $L$, we take their closures
inside a fixed algebraic closure $L^{\mathrm{alg}}$.
We write $\Gal(E):=\Gal(E^{\sep}/E)$ for the absolute
Galois group of $E$.
When $E$ has characteristic zero,
$E^{\sep}=E^{\mathrm{alg}}$; in arguments under this
assumption we use $E^{\mathrm{alg}}$, so that
$\Gal(E)=\Gal(E^{\mathrm{alg}}/E)$.

Throughout this paper, by a homomorphism between profinite
groups we mean a continuous homomorphism, and by a subgroup
of a profinite group we mean a closed subgroup.

For a profinite group $G$, let $\Imfin(G)$ denote the class
of its finite continuous quotients, up to isomorphism. We call $\Imfin(\Gal(E))$ the \emph{finite Galois-group spectrum}
of $E$; equivalently, it is the class of groups
$\Gal(L/E)$, up to isomorphism, where $L/E$ ranges over
the finite Galois extensions of $E$.

A field is \emph{pseudo-algebraically closed} (PAC) if every
nonempty geometrically integral variety over it has a rational
point. A \emph{Frobenius field} is a PAC field whose absolute
Galois group has the embedding property
\cite[Section~24.1]{FJ}.
The embedding property is defined in
Subsection~\ref{subsec:notation}.

For the fixed Melnikov formation $\D$, put
\[
\D^{\perp}=
\{H\text{ finite}\mid H\text{ has no nontrivial quotient belonging to }\D\}.
\]
This is the class denoted by $\mathcal C'$ in
\cite[Definition~1.4]{FHQuasi}, with $\mathcal C=\D$.
Let $\K_\D$ be the class of characteristic-zero Frobenius
fields $F$ such that $\Imfin(\Gal(F))=\D^{\perp}$.

For a field $E$, let $E^\D$ be the compositum, inside
$E^{\mathrm{alg}}$, of all finite Galois extensions of $E$ whose
Galois groups belong to $\D$.
For subfields $A,B$ of a common field, $AB$ denotes their
compositum. For a subfield $E$ of a field $L$, write $E_L^{\mathrm{alg}}$
for the relative algebraic closure of $E$ in $L$. A field is \emph{bounded} if, up to isomorphism over that
field, it has only finitely many finite separable extensions
of each degree.

Let $F\in\K_\D$ and let $\U$ be a monster model of $\Th(F)$.
For small elementary submodels
$M\preccurlyeq A,B\preccurlyeq\U$, define
\[
A\ind_M B
\quad\Longleftrightarrow\quad
\begin{cases}
A,B\text{ are linearly disjoint over }M,\\
(AB)_\U^{\mathrm{alg}}=(AB)^\D.
\end{cases}
\]

\begin{theorem}[Main result]\label{thm:main}
With this notation, the relation defined above is invariant
and symmetric and has existence and stationarity over models.
Put $T=\Th(F_S)$. Then:

\begin{enumerate}[label=\textup{(\roman*)},ref=(\roman*)]
\item\label{item:main-ap}
$T_\sigma$, the theory of models of $T$ with an automorphism,
has amalgamation.

\item\label{item:main-unbounded}
If $\D$ is proper, then $F$ is unbounded and
$T$ is NSOP$_1$ and not simple.

\item\label{item:main-companion}
If $\D$ is proper, then $T_\sigma$ has no model companion.
For the formation of all finite groups, the model companion
is the canonical root-predicate expansion of $\mathrm{ACFA}_0$.

\item\label{item:main-kim}
On triples of elementary submodels, this relation implies
ordinary nonforking, which implies Kim-independence, which
in turn implies linear disjointness over the base model.
If $\D$ is proper, each of these three implications is strict.

\end{enumerate}
\end{theorem}

For proper $\D$, unboundedness has a concrete Galois-theoretic
form: for every nontrivial finite quotient $H$ of $\Gal(F)$,
there are infinitely many finite Galois extensions of $F$
with Galois group isomorphic to $H$.
In particular, $\Gal(F)$ is not topologically finitely generated.
The NSOP$_1$ assertion in
Theorem~\ref{thm:main}\ref{item:main-unbounded}
is the known result for Frobenius fields
\cite[Introduction]{Chatz}, transferred to their canonical
root-predicate expansions; see Section~\ref{sec:realization}.
The non-simplicity assertion follows from our unboundedness
argument in Section~\ref{sec:unbounded}.

For example, for $\D=\{1\}$, the class $\D^{\perp}$ consists of all finite
groups and $E^\D=E$ for every field $E$.
Consequently, $\K_\D$ consists precisely of the
characteristic-zero PAC fields whose absolute Galois groups
are $\omega$-free; see Subsection~\ref{subsec:notation}.
For small elementary submodels
$M\preccurlyeq A,B\preccurlyeq\U$ as in
Theorem~\ref{thm:main}, the independence relation becomes
\[
A\ind_M B
\quad\Longleftrightarrow\quad
\begin{cases}
A,B\text{ are linearly disjoint over }M,\\
(AB)_\U^{\mathrm{alg}}=AB.
\end{cases}
\]

Our second main result concerns difference fields with a
distinguished automorphism, independently of the classes
$\K_\D$.
Such a field is \emph{pseudo-existentially closed}, or
\emph{p.e.c.}, if it is existentially closed in the ordinary
difference-field language in every difference-field extension
whose underlying field extension is regular.

\begin{theorem}[Difference fields]\label{thm:difference-main}
\begin{enumerate}[label=\textup{(\roman*)}]
\item
The class of p.e.c. characteristic-zero difference fields
is not closed under ultrapowers. In particular, it is not
elementary.
\item
The theory of
all characteristic-zero fields with an automorphism and
canonical root predicates has no model companion.
\end{enumerate}
\end{theorem}

Assertion \textup{(i)} is proved in
Theorem~\ref{thm:regular-ec-nonelementary} and
 \textup{(ii)} in
Theorem~\ref{thm:root-automorphism-no-companion}.
The root predicates make embeddings correspond to regular
extensions. Since every ordinary existential formula is also
existential in the expanded language, the reducts of its
e.c.\ models are p.e.c.\ difference fields.
The converse fails: models of $\mathrm{ACFA}_0$ are p.e.c.,
but their canonical root-predicate expansions are not e.c.\
(Remark~\ref{rem:root-ec-distinction}).

For finite $G$, we prove that a $G$-field is p.e.c. if and only
if the action is faithful and its fixed field is PAC
(Theorem~\ref{thm:finite-group-pec}).
This gives an explicit characterization of the elementary
class recorded in \cite[Remark~6.25]{Gogolok}.
A $G$-field is \emph{algebraically $G$-closed} if it has no
proper algebraic field extension to which the $G$-action extends.
For every group $G$, existential closure is equivalent to
the conjunction of p.e.c. and algebraic $G$-closedness
(Fact~\ref{prop:ec-pec-gclosed}).
We compare the axiomatizability of these three classes
for finite $G$ and for $G=\mathbb Z$.

This paper is organized as follows.
Section~\ref{sec:groups} separates notation and definitions
from the group-theoretic facts used in the proofs.
Section~\ref{sec:realization} records the algebraic-closure
and elementarity criteria, proves the regular-realization
criterion, and recalls quantifier elimination for the
canonical root-predicate expansions.
It also identifies embeddings in that language with regular
field embeddings and records the existence of existentially
closed extensions after naming an automorphism.
Sections~\ref{sec:stationarity} and~\ref{sec:automorphisms}
construct the relation, prove existence and stationarity,
and deduce amalgamation and the model-companion dichotomy.
Section~\ref{sec:unbounded} proves unboundedness and non-simplicity.
Section~\ref{sec:kim-comparison} compares the relation with
ordinary nonforking, Kim-independence, and linear disjointness,
and proves that all three implications are strict for proper
$\D$, completing the proof of Theorem~\ref{thm:main}.
Section~\ref{sec:regular-ec} begins with the closure notions
for fields with group actions and the ultraproduct criterion
for e.c.\ models.
Subsection~\ref{subsec:root-predicate-automorphisms}
proves the nonexistence of a model companion in the
root-predicate language.
Subsection~\ref{subsec:pec-difference-fields}
proves non-elementarity of p.e.c.\ difference fields,
while Subsection~\ref{sec:finite-group-actions}
gives the finite-group characterization.
Subsection~\ref{subsec:group-action-conclusions}
compares the axiomatizability of the resulting classes.

\section{Preliminaries}\label{sec:groups}

Throughout this section, $\D$ is the Melnikov formation fixed
in Section~\ref{sec:intro}.
\subsection{Notation and definitions}\label{subsec:notation}
We follow Fried--Haran \cite{FHQuasi} for formations and
residuals, and Fried--Jarden \cite{FJ} for the remaining
profinite-group and field-arithmetic terminology.

We write $N\trianglelefteq_oG$ to mean that $N$ is an open
normal subgroup of the profinite group $G$.
A profinite group is \emph{small} if it has only finitely many
open subgroups of each finite index.
For a field $E$, boundedness is equivalent to smallness of $\Gal(E)$.
A profinite group is \emph{pro-$\D$} if all its finite continuous
quotients belong to $\D$. We retain the notation $\Gal(E)$, $\Imfin(G)$, $\D^{\perp}$, and $E^\D$
introduced in Section~\ref{sec:intro}.

For a profinite group $G$, define its \emph{$\D$-residual} by
\[
\RD(G)=
\bigcap_{\substack{N\trianglelefteq_oG\\G/N\in\D}}N.
\]
\begin{remark}\label{rem:residual-spectrum}
For every profinite group $G$, one has
\[
\RD(G)=G
\quad\Longleftrightarrow\quad
\Imfin(G)\subseteq\D^{\perp}.
\]
Both conditions are equivalent to saying that $G$ has no
nontrivial finite continuous quotient belonging to $\D$.
Indeed, $\RD(G)=G$ has this meaning directly from the
definition of the residual.
If $G$ has no such quotient, then every finite quotient of
$G$ belongs to $\D^{\perp}$, since a quotient of a finite
quotient of $G$ is again a finite quotient of $G$.
Conversely, no nontrivial group in $\D$ belongs to
$\D^{\perp}$, because it is its own nontrivial $\D$-quotient.
\end{remark}
A \emph{finite embedding problem} for a profinite group $G$
is a pair of continuous epimorphisms
\[
G\xrightarrow{\alpha}H\xleftarrow{\beta}B,
\qquad H,B\text{ finite}.
\]
A \emph{weak solution} of this embedding problem is a
continuous homomorphism $\gamma:G\to B$ such that
$\beta\circ\gamma=\alpha$. It is called a \emph{proper solution}
if $\gamma$ is surjective.

A profinite group is \emph{projective} if every finite
embedding problem for it has a weak solution, and is
\emph{$\omega$-free} if every finite embedding problem
for it has a proper solution.
It has the \emph{embedding property} if every finite
embedding problem as above with $B\in\Imfin(G)$ has
a proper solution.
We write $\Fhat_\omega$ for the free profinite group of
countably infinite rank.

\subsection{Standard facts and cited results}

We first collect the residual properties used below.
They are recorded in \cite[Lemma~1.3(a),(b),(d)]{FHQuasi}.

\begin{fact}\label{fact:residual-calculus}
Let $G$ be a profinite group.
\begin{enumerate}[label=\textup{(\alph*)},ref=(\alph*)]
\item\label{item:residual-maximal}
The quotient $G/\RD(G)$ is the \emph{maximal} pro-$\D$ quotient of $G$.

\item\label{item:residual-epi}
For every epimorphism $\pi:G\twoheadrightarrow H$ of profinite
groups, we have
\[
\pi(\RD(G))=\RD(H).
\]

\item\label{item:residual-normal}
If $\RD(G)\leq V\trianglelefteq G$, then  we have
\[
\RD(V)=\RD(G).
\]

\item\label{item:residual-idempotent}
In particular,  we obtain
\[
\RD(\RD(G))=\RD(G).
\]
\end{enumerate}
\end{fact}

\begin{proof}[Proof of Fact~\ref{fact:residual-calculus}\ref{item:residual-normal}
and~\ref{item:residual-idempotent}]
For the reader's convenience, we recall the short proof.
Suppose that $\RD(G)\leq V\trianglelefteq G$.
The group $V/\RD(G)$ is pro-$\D$, since it is a closed
normal subgroup of the pro-$\D$ group $G/\RD(G)$.
By Fact~\ref{fact:residual-calculus}\ref{item:residual-maximal},
applied to $V$, we have
\[
\RD(V)\leq\RD(G).
\]

Conversely, $\RD(V)$ is characteristic in $V$ and hence
normal in $G$. The group $G/\RD(V)$ is an extension of
the pro-$\D$ group $G/V$ by the pro-$\D$ group
$V/\RD(V)$. Since $\D$ is a Melnikov formation, this
extension is again pro-$\D$.
Fact~\ref{fact:residual-calculus}\ref{item:residual-maximal},
applied to $G$, therefore gives
\[
\RD(G)\leq\RD(V).
\]
This proves part~\ref{item:residual-normal}.
Taking $V=\RD(G)$ gives part~\ref{item:residual-idempotent}.
\end{proof}

We next recall the embedding-property input.
By \cite[Lemma~4.2]{Jarden1995}, every open subgroup of an
$\omega$-free profinite group is $\omega$-free, and every
closed normal subgroup has the embedding property.
Together with the residual identities in
\cite[Lemmas~1.3 and~1.5]{FHQuasi}, this gives the following
standard consequence; compare
\cite[Theorem~5.5 and Remark~5.6(b)]{FHQuasi}.

\begin{lemma}\label{lem:group}
If $G$ is an $\omega$-free profinite group,
then $\RD(G)$ has the embedding property and
\[
\Imfin(\RD(G))=\D^{\perp}.
\]
\end{lemma}

\begin{proof}
Put $N:=\RD(G)$.
By definition, $N$ is an intersection of open normal subgroups
of $G$, so it is closed and normal in $G$.
Hence $N$ has the embedding property by
\cite[Lemma~4.2]{Jarden1995}, as recalled above.

Let $B$ be a finite quotient of $N$.
By Fact~\ref{fact:residual-calculus}\ref{item:residual-idempotent}
and~\ref{item:residual-epi}, we have $\RD(B)=B$.
Equivalently, $B$ has no nontrivial quotient in $\D$,
so $B\in\D^{\perp}$.

Conversely, suppose that $B\in\D^{\perp}$.
Since $G$ is $\omega$-free, the embedding problem
$G\to1\leftarrow B$ has a proper solution
$q:G\twoheadrightarrow B$.
Fact~\ref{fact:residual-calculus}\ref{item:residual-epi} gives
\[
q(N)=q(\RD(G))=\RD(B)=B.
\]
Thus $\Imfin(N)=\D^{\perp}$.
\end{proof}

We also record the field-arithmetic existence result used in
the regular-realization argument.
A field $E$ of characteristic zero is \emph{Hilbertian} if,
for every polynomial $f(T,X)\in E[T,X]$ that is irreducible
in $E(T)[X]$ and has positive degree in $X$, there are
infinitely many $a\in E$ such that $f(a,X)$ is irreducible
in $E[X]$ and has the same degree in $X$ as $f(T,X)$.

\begin{fact}
\label{fact:regular-pac-omega-free}
Every characteristic-zero field $E$ has a regular extension
$L/E$ such that $L$ is PAC and $\Gal(L)$ is $\omega$-free.
\end{fact}
\begin{proof}
This follows by combining the regular PAC--Hilbertian
extension construction of \cite[Proposition~13.4.6]{FJ}
with the Fried--V\"olklein theorem \cite{FV} that
characteristic-zero Hilbertian PAC fields are $\omega$-free;
see also \cite[Section~3, discussion preceding Example~3.2]{Jarden1995}.
\end{proof}

\section{Regular realization and elementary embeddings}
\label{sec:realization}

We first record the algebraic-closure and elementarity criteria
used below. We then prove the regular-realization criterion
and deduce the nonemptiness of $\K_\D$.
Finally, we record two general facts about the root-predicate
language after naming an automorphism.

For a subset $X$ of a field $F$, we write $\acl_F(X)$ for its
model-theoretic algebraic closure in $F$.

We use the root-predicate reduct of the Galois language of
Cherlin--van den Dries--Macintyre
\cite[Section~4.2, p.~60]{CDM}. Let us set
\[
\mathcal L_S
:=\mathcal L_{\mathrm{rings}}\cup\{S_n:n\geq2\},
\]
where $S_n$ is an $n$-ary predicate symbol. For a field $K$, write $K_S$ for
its canonical expansion, defined by
\[
S_n(c_1,\ldots,c_n)
\quad\Longleftrightarrow\quad
\exists x\,(x^n+c_1x^{n-1}+\cdots+c_n=0).
\]
Thus $K_S$ is a parameter-free definitional expansion of $K$.
In particular, the expansion does not change the definable
sets or model-theoretic algebraic closure.
Every field automorphism preserves this expansion.

\begin{proposition}\label{prop:criteria}
The following assertions hold in characteristic zero.
\begin{enumerate}[label=\textup{(\alph*)}]
\item If $F$ is PAC and $A\subseteq F$ is a subfield, then
$\acl_F(A)=A_F^{\mathrm{alg}}$.
\item If $F\subseteq K$ are Frobenius fields, $K/F$ is regular,
and $\Imfin(\Gal(F))=\Imfin(\Gal(K))$, then $F\preccurlyeq K$.
\item If $F$ is a Frobenius field and $A,B$ are relatively
algebraically closed subfields of $F$, then every field isomorphism
$f:A\to B$ is a partial elementary map of $F$.
\end{enumerate}
The same elementarity conclusions hold for the canonical
root-predicate expansions.
\end{proposition}

\begin{proof}
Assertion~(a) is the algebraic-closure description recalled in
\cite[Section~1.8]{Chatz}.

For~(b) and~(c), we recall the Galois formalism of
Cherlin--van den Dries--Macintyre
\cite[Section~4.2, p.~60]{CDM}:
\[
\mathcal L_{\mathrm{Gal}}
=
\mathcal L_S\cup
\{I_G:G\text{ a finite group up to isomorphism}\}.
\]
The predicates $S_n$ have the root interpretations given
above, and $I_G$ is a nullary predicate interpreted by
\[
K\models I_G
\quad\Longleftrightarrow\quad
G\in\Imfin(\Gal(K)).
\]
In the terminology of \cite{CDM}, perfect Iwasawa RC-fields
are precisely perfect Frobenius fields.
Their theory $\mathrm{pIRC}$, with these interpretations,
eliminates quantifiers in $\mathcal L_{\mathrm{Gal}}$
by \cite[Theorem~41, p.~60]{CDM}.
This theorem does not require a fixed characteristic or
finite Galois-group spectrum.
In particular, it applies to characteristic-zero
Frobenius fields.

For~(b), regularity implies that the inclusion $F\subseteq K$
preserves and reflects every root predicate.
The equality of the finite Galois-group spectra ensures
that the nullary predicates agree.
Thus the inclusion is an $\mathcal L_{\mathrm{Gal}}$-embedding,
and quantifier elimination makes it elementary.

For~(c), a polynomial over $A$ has a root in $F$ if and only
if it has a root in $A$, because $A$ is relatively
algebraically closed in $F$. The same holds for $B$.
Therefore $f$ preserves and reflects the root predicates
inherited from $F$.
The substructures on $A$ and $B$ also inherit the same
nullary predicates from $F$.
Thus $f$ is a partial $\mathcal L_{\mathrm{Gal}}$-isomorphism,
and quantifier elimination makes it partial elementary.
Taking reducts gives the conclusions in both the ring
language and $\mathcal L_S$.
\end{proof}

By Proposition~\ref{prop:criteria}(a), every subfield
$E\subseteq\U$ satisfies
\[
E_\U^{\mathrm{alg}}=\acl_\U(E),
\]
since $\U$ is PAC.
Thus the closure condition in the introductory definition
of $\ind$ can equivalently be written as
\[
\acl_\U(AB)=(AB)^\D.
\]
We use this equivalent formulation in the model-theoretic
arguments below.
The same identification holds in the canonical root-predicate
expansion, which is definitional.

We next record the Galois-theoretic interpretation of the
field construction $E^\D$.

\begin{fact}
\label{fact:galois-residual}
Let $E$ be a characteristic-zero field.
\begin{enumerate}[label=\textup{(\alph*)},ref=(\alph*)]
\item\label{item:galois-fixed-field}
The extension $E^\D/E$ is Galois, and
\[
E^\D=(E^{\mathrm{alg}})^{\RD(\Gal(E))}.
\]
The restriction map induces the following canonical isomorphism
\[
\Gal(E)/\RD(\Gal(E))\cong\Gal(E^\D/E).
\]

\item\label{item:galois-absolute-group}
We have
\[
\Gal(E^\D)=\RD(\Gal(E)),
\qquad
E^\D=E
\quad\Longleftrightarrow\quad
\RD(\Gal(E))=\Gal(E).
\]
\end{enumerate}
\end{fact}

\begin{proof}
By definition, $E^\D$ is the compositum in $E^{\mathrm{alg}}$ of
the finite Galois extensions $L/E$ with $\Gal(L/E)\in\D$.
An element of $\Gal(E)$ fixes this compositum if and only if it
fixes every such $L$. Therefore
\[
\Gal(E^{\mathrm{alg}}/E^\D)
=
\bigcap_{\substack{N\trianglelefteq_o\Gal(E)\\\Gal(E)/N\in\D}}N
=
\RD(\Gal(E)).
\]
Galois correspondence gives
part~\ref{item:galois-fixed-field}.
Since $E^\D/E$ is algebraic, $E^{\mathrm{alg}}$
is also an algebraic closure of $E^\D$, proving the equality
in part~\ref{item:galois-absolute-group}.
\end{proof}

The same description gives the following compatibility
with regular extensions.

\begin{lemma}
\label{lem:residual-regular-intersection}
If $E'/E$ is a regular extension of characteristic-zero
fields, then we have
\[
(E')^\D\cap E^{\mathrm{alg}}=E^\D.
\]
\end{lemma}

\begin{proof}
Regularity makes the restriction map
\[
\rho:\Gal(E')\longrightarrow\Gal(E)
\]
surjective.
By Fact~\ref{fact:residual-calculus}\ref{item:residual-epi},
we obtain
\[
\rho\bigl(\RD(\Gal(E'))\bigr)=\RD(\Gal(E)).
\]
Using the fixed-field descriptions from
Fact~\ref{fact:galois-residual}\ref{item:galois-fixed-field}
for both $E'$ and $E$, we obtain
\[
(E')^\D\cap E^{\mathrm{alg}}
=
(E^{\mathrm{alg}})^{\rho(\RD(\Gal(E')))}
=
(E^{\mathrm{alg}})^{\RD(\Gal(E))}
=
E^\D. \qedhere
\]
\end{proof}

The following regular-realization criterion follows from
Fact~\ref{fact:regular-pac-omega-free},
Lemma~\ref{lem:group}, and the residual identities in
Fact~\ref{fact:residual-calculus}.
We include the short argument for completeness.

\begin{lemma}\label{lem:realize}
Let $E$ be a characteristic-zero field. The following are equivalent:
\begin{enumerate}[label=\textup{(\roman*)}]
\item $E$ has a regular extension $K\in\K_\D$;
\item $\RD(\Gal(E))=\Gal(E)$;
\item $\Imfin(\Gal(E))\subseteq\D^{\perp}$.
\end{enumerate}
\end{lemma}

\begin{proof}
The equivalence of \textup{(ii)} and \textup{(iii)}
follows from Remark~\ref{rem:residual-spectrum}.
 If (i) holds, restriction
$\Gal(K)\twoheadrightarrow \Gal(E)$ is surjective, so every finite
quotient of $\Gal(E)$ belongs to
$\Imfin(\Gal(K))=\D^{\perp}$. Thus (i) implies (iii).

Assume \textup{(ii)}.
By Fact~\ref{fact:regular-pac-omega-free}, choose a regular
PAC extension $L/E$ with $\Gal(L)$ $\omega$-free.
Put $K=L^\D$.
As an algebraic extension of a PAC field, $K$ is PAC
\cite[Corollary~11.2.5]{FJ}.
By Fact~\ref{fact:galois-residual}\ref{item:galois-absolute-group},
\[
\Gal(K)=\RD(\Gal(L)).
\]
Lemma~\ref{lem:group} therefore gives the embedding property
and $\Imfin(\Gal(K))=\D^{\perp}$, hence $K\in\K_\D$.

Finally, assumption~\textup{(ii)} and
Fact~\ref{fact:galois-residual}\ref{item:galois-absolute-group}
give $E^\D=E$.
Since $L/E$ is regular,
Lemma~\ref{lem:residual-regular-intersection} gives
\[
K\cap E^{\mathrm{alg}}
=L^\D\cap E^{\mathrm{alg}}
=E^\D=E.
\]
Hence $K/E$ is regular.
\end{proof}
Taking $E=\mathbb Q^{\operatorname{alg}}$, we have
\[
\Gal(E)=1=\RD(\Gal(E)).
\]
Lemma~\ref{lem:realize} therefore supplies a regular extension
$K/E$ belonging to $\K_\D$. In particular, $\K_\D$ is nonempty.

\begin{remark}
By Fact~\ref{fact:galois-residual}\ref{item:galois-absolute-group},
condition \textup{(ii)} of Lemma~\ref{lem:realize} is equivalent
to $E^\D=E$.
Thus a characteristic-zero field $E$ satisfies $E^\D=E$ if and
only if it occurs as a relatively algebraically closed subfield
of some $K\in\K_\D$.
By Proposition~\ref{prop:criteria}(a), this is equivalently the
condition $\acl_K(E)=E$, where $\acl_K$ is model-theoretic
algebraic closure.
The same characterization holds in the canonical root-predicate
expansion of $K$, since that expansion is definitional.
\end{remark}

We now axiomatize the classes of fields under consideration
and spell out the quantifier-elimination result used above.

For each finite group $G$, let $\varphi_G$ be a sentence
in the ring language expressing that the field has a finite
Galois extension with Galois group isomorphic to $G$.
The requirement $\Imfin(\Gal(K))=\D^{\perp}$ is expressed by the
axiom scheme
\[
\Sigma_{\D}
=\{\varphi_G:G\in\D^{\perp}\}
 \cup\{\neg\varphi_G:G\notin\D^{\perp}\},
\]
where $G$ ranges over finite groups up to isomorphism.
Let $T_\D$ be the theory in the ring language consisting
of the axioms for characteristic-zero Frobenius fields
\cite[Section~1.9]{Chatz} together with $\Sigma_\D$.
Then
\[
\operatorname{Mod}(T_\D)=\K_\D.
\]
In particular, $\K_\D$ is an elementary class.
Let $T_{\D,S}$ denote the definitional expansion of $T_\D$
by the root predicates defined above. Its models are
precisely the structures $K_S$ with $K\in\K_\D$.

\begin{fact}
\label{fact:root-qe}
Let $\mathcal C$ be a class of finite groups, up to isomorphism,
realized as the finite Galois-group spectrum of a
characteristic-zero Frobenius field.
The theory of characteristic-zero Frobenius fields with
spectrum $\mathcal C$, in their canonical $\mathcal L_S$-expansions,
is consistent and eliminates quantifiers.
In particular, $T_{\D,S}$ is consistent, eliminates quantifiers,
and is model complete.
\end{fact}

\begin{proof}
Use the language $\mathcal L_{\mathrm{Gal}}$ and the
quantifier-elimination theorem recalled in the proof of
Proposition~\ref{prop:criteria}.
Over the canonical expansions of characteristic-zero
Frobenius fields, every $\mathcal L_S$-formula is equivalent
to a quantifier-free $\mathcal L_{\mathrm{Gal}}$-formula.
When the finite Galois-group spectrum is $\mathcal C$,
replace each nullary predicate $I_G$ in this formula by
$0=0$ if $G\in\mathcal C$, and by $0\ne0$ otherwise.
The resulting formula is quantifier-free in $\mathcal L_S$.
Consistency follows from the assumed realization of
$\mathcal C$.

Lemma~\ref{lem:realize}, applied to
$E=\mathbb Q^{\operatorname{alg}}$, shows that
$\D^\perp$ is realized by a field in $\K_\D$.
The first assertion therefore applies to $T_{\D,S}$.
Model completeness follows from quantifier elimination.
\end{proof}
By Fact~\ref{fact:root-qe}, every $\mathcal L_S$-embedding
between models of $T_{\D,S}$ is elementary.
Since regular field embeddings preserve and reflect the
root predicates, every regular embedding between fields
in $\K_\D$ is elementary, both in $\mathcal L_S$ and in
the ring language.

The theory $T_{\D,S}$ need not be complete.
For a fixed $F\in\K_\D$, we work with its completion
\[
T=\Th(F_S).
\]
This is a parameter-free definitional expansion of $\Th(F)$.
It inherits quantifier elimination and model completeness
from $T_{\D,S}$.

The theory of every Frobenius field is NSOP$_1$; see
\cite[Introduction and Comments~2.13]{Chatz}, which recalls
the relevant results of Chernikov--Ramsey.
Consequently, for every $F\in\K_\D$, the theory
$T=\Th(F_S)$ is NSOP$_1$, since definitional expansions
preserve this property.
This proves the NSOP$_1$ assertion in
Theorem~\ref{thm:main}\ref{item:main-unbounded}.
We also use this result in the comparison with
Kim-independence in Section~\ref{sec:kim-comparison}.

We finish with two facts about fields with an automorphism
in the root-predicate language.

For later use, put
\[
\mathcal L_{S,\sigma}=\mathcal L_S\cup\{\sigma\},
\]
where $\sigma$ is a unary function symbol.
Let $T_{S,\sigma}$ be the theory of all characteristic-zero
fields with an automorphism and the canonical interpretations
of the root predicates.
We write $(K,\sigma)_S$ for the canonical expansion of a
difference field $(K,\sigma)$ to this language.
No PAC condition or prescribed finite Galois-group spectrum
is imposed on models of $T_{S,\sigma}$.

\begin{proposition}\label{prop:root-regular-embeddings}
For an inclusion $K\subseteq L$ of characteristic-zero fields,
\[
\begin{aligned}
K_S\subseteq_{\mathcal L_S}L_S
&\quad\Longleftrightarrow\quad
K\text{ is relatively algebraically closed in }L\\
&\quad\Longleftrightarrow\quad
L/K\text{ is regular}.
\end{aligned}
\]
Consequently, embeddings between models of $T_{S,\sigma}$ are
exactly the difference-field embeddings whose underlying field
extensions are regular.
\end{proposition}

\begin{proof}
If $a\in L\setminus K$ is algebraic over $K$, its monic minimal
polynomial has degree $n\geqslant2$. Its coefficient tuple satisfies
$S_n$ in $L$ but not in $K$, contradicting substructure inclusion.
Conversely, if $K$ is relatively algebraically closed in $L$,
every root in $L$ of a polynomial over $K$ already lies in $K$.
Thus all root predicates are preserved and reflected.
In characteristic zero this condition is equivalent to regularity.
The difference-field assertion follows by additionally requiring
the embedding to commute with the automorphisms.
\end{proof}

\begin{fact}
\label{lem:root-ec-hull}
The theory $T_{S,\sigma}$ is inductive.
Every model $(K,\sigma)_S$ has an existentially closed
$T_{S,\sigma}$-extension $(M,\tau)_S$.
The underlying field extension $M/K$ is regular.
\end{fact}

\begin{proof}
The field, characteristic-zero, automorphism, and root-predicate
axioms admit a universal--existential axiomatization.
Hence $T_{S,\sigma}$ is inductive.
Every model of an inductive theory has an existentially closed
extension. The underlying field extension is regular by
Proposition~\ref{prop:root-regular-embeddings}.
\end{proof}

\section{Existence and stationarity over models}\label{sec:stationarity}

In this section, we prove the assertion about the relation
in Theorem~\ref{thm:main}.
Recall that $\D$ is a Melnikov formation: a nonempty class
of finite groups closed under isomorphisms, normal subgroups,
quotients, and extensions.
Fix $F\in\K_\D$; thus $F$ is a characteristic-zero Frobenius
field with $\Imfin(\Gal(F))=\D^\perp$, where $\D^\perp$
consists of the finite groups having no nontrivial quotient
in $\D$.
We work in a sufficiently saturated and homogeneous
monster $\U\models\Th(F)$.
For small elementary submodels
$M\preccurlyeq A,B\preccurlyeq\U$, recall that
\[
A\ind_M B
\quad\Longleftrightarrow\quad
\begin{cases}
A,B\text{ are linearly disjoint over }M,\\
\acl_\U(AB)=(AB)^\D.
\end{cases}
\]
\begin{proposition}[Existence]\label{prop:exist}
For small models $M\preccurlyeq A,B\preccurlyeq\U$, there is $A'\preccurlyeq\U$ such that
\[
A'\cong_M A,\qquad A'\ind_M B.
\]
\end{proposition}

\begin{proof}
Since $A/M$ and $B/M$ are regular, $A\otimes_M B$ is a domain.
Put
\[
E=\Frac(A\otimes_M B).
\]
We identify $A$ and $B$ with their canonical images in $E$.
These images are linearly disjoint over $M$, and the
extensions $E/A$ and $E/B$ are regular.
By Fact~\ref{fact:galois-residual}\ref{item:galois-absolute-group},
\[
\Gal(E^\D)=\RD(\Gal(E)).
\]
Since $A,B\in\K_\D$, we have $A^\D=A$ and $B^\D=B$.
Applying Lemma~\ref{lem:residual-regular-intersection}
to the regular extensions $E/A$ and $E/B$, we obtain
\[
E^\D\cap A^{\mathrm{alg}}=A^\D=A,
\qquad
E^\D\cap B^{\mathrm{alg}}=B^\D=B.
\]
Thus $E^\D/A$ and $E^\D/B$ are regular. Moreover,
Fact~\ref{fact:residual-calculus}\ref{item:residual-idempotent}
gives
\[
\RD(\Gal(E^\D))
=\RD(\RD(\Gal(E)))
=\RD(\Gal(E))
=\Gal(E^\D).
\]
We apply Lemma~\ref{lem:realize} to obtain a regular extension $K/E^\D$ such that $K\in\K_\D$.
The inclusions $A,B\hookrightarrow K$ are regular, hence
elementary by Proposition~\ref{prop:criteria}(b).
As $E^\D/E$ is algebraic and $K/E^\D$ is regular,
Proposition~\ref{prop:criteria}(a) gives
\[
\acl_K(E)=E^\D.
\]
Replacing $K$ by an elementary submodel containing $E^\D$,
we may assume that $K$ is small; the displayed equality
and regularity over $E^\D$ are preserved.
Choose an elementary embedding $j:K\to\U$ over $B$, and
put $A'=j(A)$ and $f=j|_A$.
Then $A'\preccurlyeq\U$, and $f:A\to A'$ is an elementary
isomorphism fixing $M$ pointwise.
The fields $A'$ and $B$ are linearly disjoint over $M$.
Since algebraic closure of a subset of an elementary
submodel is unchanged in the monster, we have
\[
\acl_\U(A'B)=j(E^\D)=(A'B)^\D.
\]
Thus $A'\ind_M B$, as required.
\end{proof}

\begin{proposition}[Stationarity]\label{prop:stat}
Suppose $A\ind_M B$ and $A'\ind_M B$. If $f:A\to A'$ is an elementary isomorphism over $M$, then
\[
f\cup\mathrm{id}_B:A\cup B\to A'\cup B
\]
 is a partial elementary map.
\end{proposition}

\begin{proof}
Linear disjointness extends $f\cup\mathrm{id}_B$ uniquely to a field isomorphism
\[
f_0:E:=AB\longrightarrow E':=A'B.
\]
We extend $f_0$ further to an isomorphism
$\widehat f:E^{\mathrm{alg}}\to (E')^{\mathrm{alg}}$.
Conjugation by $\widehat f$ induces an isomorphism
$\Gal(E)\to \Gal(E')$, which sends $\RD(\Gal(E))$ onto $\RD(\Gal(E'))$
by Fact~\ref{fact:residual-calculus}\ref{item:residual-epi}.
Taking fixed fields and applying
Fact~\ref{fact:galois-residual}\ref{item:galois-fixed-field},
we obtain
\[
\widehat f(E^\D)=(E')^\D.
\]
Since $A\ind_M B$ and $A'\ind_M B$, we obtain
\[
E^\D= \acl_\U(E)=\acl_\U(AB),\qquad (E')^\D= \acl_\U(E')=\acl_\U(AB).
\]
The map $\widehat f|_{E^\D}$ is therefore an isomorphism between relatively algebraically closed subfields of $\U$, so it is partial elementary by Proposition~\ref{prop:criteria}(c). Since we have
\[
\left(\widehat f|_{E^\D}\right)|_{A\cup B}=f\cup\mathrm{id}_B,
\]
the map $f\cup\mathrm{id}_B$ is partial elementary as well.
\end{proof}

Invariance and symmetry follow directly from the definition
of the relation.
This proves the invariance, symmetry, existence and stationarity
asserted in Theorem~\ref{thm:main}.

\section{No model companion after adding an automorphism}\label{sec:automorphisms}

For a complete theory $T$ in a language $L$, let $T_\sigma$
be the theory of models of $T$ equipped with an $L$-automorphism.
We say that $T$ has the \emph{amalgamation property for
automorphisms}, abbreviated \emph{PAPA}, if such structures
have amalgamation with respect to maps that are elementary
on the underlying $L$-structures and commute with the named
automorphisms; see \cite[Definition~3.1]{KikyoPAPA}.
If $T$ is model complete, PAPA is equivalent to the ordinary
amalgamation property of $T_\sigma$.

The following criterion isolates the independence properties
used in the standard proof of PAPA for stable theories;
compare \cite[Fact~3.5]{KikyoPAPA}.

\begin{lemma}
\label{lem:papa-criterion}
Let $T$ be a complete theory with monster model $\U$.
Suppose that a relation $A\indgen_M B$ is defined for small
elementary submodels $M\preccurlyeq A,B\preccurlyeq\U$
and satisfies the following properties.
\begin{enumerate}[label=\textup{(\roman*)}]
\item
\emph{Invariance.}
For every $g\in\operatorname{Aut}(\U)$,
\[
A\indgen_M B
\quad\Longrightarrow\quad
g(A)\indgen_{g(M)}g(B).
\]

\item
\emph{Existence.}
For every $M\preccurlyeq A,B\preccurlyeq\U$, there is
an elementary copy $A'$ of $A$ over $M$ such that
$A'\indgen_M B$.

\item
\emph{Stationarity.}
If $A\indgen_M B$ and $A'\indgen_M B$, then every
elementary isomorphism $f:A\to A'$ fixing $M$ pointwise
has the property that
$f\cup\operatorname{id}_B$ is a partial elementary map.
\end{enumerate}
Then $T$ has PAPA.
\end{lemma}

\begin{proof}
We consider models $(M,\sigma_M)$, $(A,\sigma_A)$ and
$(B,\sigma_B)$ of $T_\sigma$ such that
\[
M\preccurlyeq A,B,
\qquad
\sigma_A|_M=\sigma_B|_M=\sigma_M;
\]
and we work in a sufficiently large monster model $\U$ of $T$.
By existence, we may replace $A$ by an elementary copy
over $M$, transporting $\sigma_A$, so that $A\indgen_M B$.

We extend $\sigma_B$ to an automorphism
$\tau\in\operatorname{Aut}(\U)$, and define
\[
A':=\tau^{-1}(A),
\qquad
f:=\tau^{-1}\circ\sigma_A:A\longrightarrow A'.
\]
Since $\tau(M)=M$ and $\tau(B)=B$, invariance gives
$A'\indgen_M B$.
Moreover, $f$ is an elementary isomorphism fixing $M$
pointwise, because
$\tau|_M=\sigma_B|_M=\sigma_A|_M$.
Stationarity therefore implies that
$f\cup\operatorname{id}_B$ is partial elementary.
Composing with $\tau$, we obtain the partial elementary map
\[
\tau\circ(f\cup\operatorname{id}_B)
=
\sigma_A\cup\sigma_B.
\]
By homogeneity, this map extends to an automorphism $s$
of $\U$. Thus $(\U,s)$ is an amalgam.
\end{proof}

The criterion requires only invariance, existence and
stationarity over models.
Symmetry, monotonicity, transitivity, finite character
and local character are not used.
In particular, the relation need not be identified with
forking or Kim-independence.
For a complete stable theory, nonforking over models
satisfies the three hypotheses, recovering stable PAPA.

We now return to the fields considered in this paper.
Let $T=\Th(F_S)$ for $F\in\K_\D$, with the canonical
root-predicate expansion defined in Section~\ref{sec:realization}.
By quantifier elimination, $T$ is model complete.
Its embeddings between models are exactly the regular
field embeddings, and these are elementary on the field
reducts.

\begin{proposition}\label{prop:ap}
The theory $T$ has PAPA.
Equivalently, $T_\sigma$ has the amalgamation property.
\end{proposition}

\begin{proof}
The relation $\ind$ is invariant and has existence and
stationarity over models by Propositions~\ref{prop:exist} and~\ref{prop:stat}.
These properties also hold in the canonical
root-predicate expansion, since it is definitional.
Lemma~\ref{lem:papa-criterion} therefore gives PAPA,
and model completeness gives the amalgamation property
of $T_\sigma$.
\end{proof}

\begin{corollary}\label{cor:noMC}
If $\D$ is proper, then $T$ has the independence property
and $T_\sigma$ has no model companion.
\end{corollary}

\begin{proof}
A proper Melnikov formation omits some finite simple group $S$.
Such an $S$ belongs to $\D^{\perp}$, so $\Gal(F)$ is nontrivial and
$F$ is not algebraically closed.
A PAC field of characteristic zero which is not algebraically
closed has the independence property
\cite[Introduction]{Chatz}.
This property is preserved by the canonical root-predicate
expansion, which is definitional.
Thus $T$ is unstable.

The theory $T$ is model complete by
Fact~\ref{fact:root-qe}, and $T_\sigma$ has amalgamation
by Proposition~\ref{prop:ap}.
Kikyo's theorem states that if a complete theory is unstable
and model complete and its automorphism theory has
amalgamation, then that automorphism theory has no model
companion \cite[Theorem~3.3]{Kikyo}.
Applying this theorem to $T$ proves the conclusion.
\end{proof}

\section{Unboundedness and non-simplicity}\label{sec:unbounded}

In this section, we prove the unboundedness and non-simplicity assertions in
Theorem~\ref{thm:main}.

\begin{lemma}\label{lem:powers}
Let $G$ be a profinite group, and let $H$ be a nontrivial
finite group. Suppose that $H^n$ is a continuous quotient
of $G$ for every positive integer $n$.
Then there are infinitely many open normal subgroups
$N\trianglelefteq_o G$ such that $G/N\cong H$.
In particular, $G$ is not small.
\end{lemma}

\begin{proof}
For every $n\geqslant1$, choose an epimorphism
$q_n:G\twoheadrightarrow H^n$.
Composing with the $n$ coordinate projections gives
epimorphisms from $G$ onto $H$ with pairwise distinct kernels.
Thus, for every $n$, there are at least $n$ open normal
subgroups $N$ with $G/N\cong H$.
\end{proof}

\begin{theorem}\label{thm:infinite}
Let $G$ be a nontrivial profinite group such that
\[
\Imfin(G)=\D^{\perp}
\]
for a Melnikov formation $\D$.
Then, for every nontrivial finite quotient $H$ of $G$, there are
infinitely many open normal subgroups $N\trianglelefteq_o G$
such that $G/N\cong H$.
In particular, $G$ is not small.
\end{theorem}

\begin{proof}
Let us fix a nontrivial finite quotient $H$ of $G$. By Lemma~\ref{lem:powers}, it is enough to show that $H^n\in\D^{\perp}$ for every positive integer $n$.
Suppose that
\[
q:H^n\twoheadrightarrow Q,\qquad Q\in\D.
\]
We aim to show that $Q$ is trivial. For each coordinate subgroup $H_i$ of $H^n$, its image
$q(H_i)$ is normal in $Q$.
A Melnikov formation is closed under normal subgroups,
so $q(H_i)\in\D$.
On the other hand, $q(H_i)$ is a quotient of $H$, and
$H\in\D^{\perp}$. Hence $q(H_i)=1$ for every $i$.
The coordinate subgroups generate $H^n$, so $Q=1$ as desired.
\end{proof}

\begin{corollary}
\label{cor:unbounded}
If $\D$ is proper and $F\in\K_\D$, then $\Gal(F)$ is not small and $\Th(F)$ is not simple.
\end{corollary}

\begin{proof}
By the definition of the class $\K_\D$, we have
\[
\Imfin(\Gal(F))=\D^{\perp}.
\]
Since $\D$ is proper, $\Gal(F)$ is non-trivial and Theorem \ref{thm:infinite} implies that $\Gal(F)$ is not small.
In particular, $F$ is unbounded. Every simple PAC field is bounded
\cite[Theorem~3.9]{ChatzSimple}, so
 $\Th(F)$ is not simple. The same holds for $T=\Th(F_S)$, since the root-predicate
expansion is definitional.
\end{proof}

Corollary~\ref{cor:unbounded}, together with the NSOP$_1$
result recalled in Section~\ref{sec:realization}, proves
part~\ref{item:main-unbounded} of Theorem~\ref{thm:main}.
Sections~\ref{sec:stationarity} and~\ref{sec:automorphisms}
establish parts~\ref{item:main-ap} and~\ref{item:main-companion},
as well as the assertions concerning invariance, symmetry,
existence and stationarity.
Part~\ref{item:main-kim} is proved in
Section~\ref{sec:kim-comparison}.

\section{Comparison of $\D$-independence with forking and Kim-independence}
\label{sec:kim-comparison}

Let $T$ be the complete theory of the canonical root-predicate
expansion of a field $F\in\K_\D$, and work in a monster model
$\U\models T$.
Write $\indfork$ for ordinary nonforking and $\indKim$
for Kim-independence.
Since the root-predicate expansion is definitional, these
relations agree with their counterparts in the field reduct.
For subfields $M\subseteq A,B$ of a common field, write
$A\indld_M B$ if $A$ and $B$ are linearly disjoint over $M$. We often write $\operatorname{tp}(A/M)$ implicitly meaning some enumeration of $A$.

For elementary submodels $M\preccurlyeq A,B\preccurlyeq\U$,
the definition of our relation becomes
\[
A\ind_M B
\quad\Longleftrightarrow\quad
A\indld_M B
\ \text{ and }\
\acl_\U(AB)=(AB)^\D.
\]
We obtain the implication from this relation to ordinary
nonforking by constructing global invariant extensions.
Lemma~\ref{lem:stationary-invariant-extension} records the
abstract argument for relations on elementary submodels,
and Proposition~\ref{prop:residual-invariant-extension}
applies it to our relation.
Kim-independence implies linear disjointness and, in our
characteristic-zero setting, a further condition on the
algebraic closures, recalled in
Fact~\ref{fact:kim-linear-disjointness}.
Among Kim-independent triples of elementary models,
our relation is characterized by the additional equality
$\acl_\U(AB)=(AB)^\D$.
Theorem~\ref{thm:residual-kim-comparison} proves that every
implication in this comparison chain is strict when $\D$
is proper, whereas all four conditions coincide in the
algebraically closed case.

\begin{lemma}
\label{lem:residual-model-monotonicity}
If $M\preccurlyeq B\preccurlyeq C\preccurlyeq\U$,
$M\preccurlyeq A\preccurlyeq\U$, and $A\ind_M C$, then $A\ind_M B$.
\end{lemma}

\begin{proof}
The assumption gives $A\indld_M C$, hence $A\indld_M B$.

Moreover, $AC/AB$ is regular. Indeed, $C/B$ is regular because
$B\preccurlyeq C$. Linear disjointness identifies
$AB\otimes_B C$ with a localization of $A\otimes_M C$ and
identifies its fraction field with $AC$. Thus $AC/AB$ is
obtained by regular base change from $C/B$.

Using Proposition~\ref{prop:criteria}(a), the assumption
$A\ind_M C$, and
Lemma~\ref{lem:residual-regular-intersection}, we obtain
\[
\begin{aligned}
\acl_\U(AB)
&=\acl_\U(AC)\cap(AB)^{\mathrm{alg}}\\
&=(AC)^\D\cap(AB)^{\mathrm{alg}}\\
&=(AB)^\D.
\end{aligned}
\]
Together with $A\indld_M B$, this proves $A\ind_M B$.
\end{proof}

The next lemma is the version for elementary submodels of the
stationary-extension construction in
\cite[Fact~3.2(1) in the arXiv version]{MutchnikConant}.
We include the proof to specify the hypotheses needed here.

\begin{lemma}
\label{lem:stationary-invariant-extension}
Let $T_0$ be a complete theory with monster model $\U$.
Suppose that a relation $A\indgen_M B$, defined for small
elementary submodels $M\preccurlyeq A,B\preccurlyeq\U$,
satisfies invariance, existence and stationarity over models
in the sense of Lemma~\ref{lem:papa-criterion}.
Assume in addition \emph{right monotonicity}:
\[\text{for all\
$M\preccurlyeq N\preccurlyeq N'\preccurlyeq\U$}\qquad \qquad
A\indgen_M N'
\quad\Longrightarrow\quad
A\indgen_M N.
\]
If $A\indgen_M B$, then the type $\operatorname{tp}(A/B)$ has a
global $M$-invariant extension and
$A\indfork_MB$.
\end{lemma}

\begin{proof}
For every small elementary model
$M\preccurlyeq N\preccurlyeq\U$,
existence supplies
a model $A_N\cong_MA$ (the isomorphism is a partial elementary map) such that $A_N\indgen_M N$.
Stationarity implies that
\[
p_N:=\operatorname{tp}(A_N/N)
\]
does not depend on the choice of $A_N$.

If $M\preccurlyeq N\preccurlyeq N'\preccurlyeq\U$,
right monotonicity gives $A_{N'}\indgen_M N$.
Stationarity therefore yields
\[
p_{N'}|_N=p_N.
\]
Thus the types $p_N$ form a coherent family, and their union
\[
q:=
\bigcup_{\substack{M\preccurlyeq N\preccurlyeq\U\\
                  N\text{ small}}}p_N
\]
is consistent. Every finite tuple of parameters from
$\U$ is in one such $N$, so $q$ is complete
over $\U$.

For $\tau\in\operatorname{Aut}(\U/M)$,
invariance and stationarity give
\[
\tau(p_N)=p_{\tau(N)}.
\]
Consequently $\tau(q)=q$, and $q$ is $M$-invariant.
Its restrictions have the stated independence property
by construction and invariance.
Any other global extension with this property must restrict
to $p_N$ over every $N$, so it equals $q$.
Since $A\indgen_M B$, stationarity gives
\[
q|_B=\operatorname{tp}(A/B),
\]
A type admitting a global $M$-invariant extension
does not fork over $M$; see \cite[Fact~7.2(3)]{KRKim}.
\end{proof}

\begin{proposition}
\label{prop:residual-invariant-extension}
Let $M\preccurlyeq A,B\preccurlyeq\U$ and suppose that
$A\ind_M B$. Then, for every small tuple $\bar a$ from $A$,
the type $\operatorname{tp}(\bar a/B)$ extends to a global
$M$-invariant type. Consequently,
\[
A\ind_M B
\quad\Longrightarrow\quad
A\mathrel{\indfork_M}B
\quad\Longrightarrow\quad
A\mathrel{\indKim_M}B.
\]
\end{proposition}

\begin{proof}
Apply Lemma~\ref{lem:stationary-invariant-extension}
to the relation $\ind$.
Invariance follows from its definition.
Existence and stationarity are
Propositions~\ref{prop:exist} and~\ref{prop:stat},
and right monotonicity is
Lemma~\ref{lem:residual-model-monotonicity}.
The lemma gives the required invariant extensions and
ordinary nonforking.
Ordinary nonforking implies Kim-independence over models
by \cite[Fact~7.2(3)]{KRKim}.
\end{proof}

\begin{fact}
\label{fact:kim-linear-disjointness}
In any complete theory expanding the theory of fields,
let $M\preccurlyeq A,B$ be elementary submodels of a common
monster model $\U$.
If $A\indKim_M B$, then $A\indld_M B$ and
\[
\acl_{\U}(AB)\cap A^{\sep}B^{\sep}=AB.
\]
\end{fact}

\begin{proof}
Both assertions follow from
\cite[Proposition~9.28(1),(3)]{KRKim}.
\end{proof}

\begin{corollary}
\label{cor:residual-kim-closure}
For the theory $T=\Th(F_S)$ considered in this section,
let $M\preccurlyeq A,B\preccurlyeq\U$.
The following conditions are equivalent:
\begin{enumerate}[label=\textup{(\roman*)}]
\item
$A\ind_M B$;
\item
$A\indfork_M B$ and $\acl_\U(AB)=(AB)^\D$;
\item
$A\indKim_M B$ and $\acl_\U(AB)=(AB)^\D$.
\end{enumerate}
\end{corollary}

\begin{proof}
Proposition~\ref{prop:residual-invariant-extension} and the
definition of $\ind$ give \textup{(i)}$\Rightarrow$\textup{(ii)}.
Ordinary nonforking implies Kim-independence, giving
\textup{(ii)}$\Rightarrow$\textup{(iii)}.
If \textup{(iii)} holds, Fact~\ref{fact:kim-linear-disjointness}
gives $A\indld_M B$.
Together with the assumed closure equality, this is precisely
the definition of \textup{(i)}.
\end{proof}

Consequently, on triples
$M\preccurlyeq A,B\preccurlyeq\U$, the relations satisfy
\[
\ind
\quad\Longrightarrow\quad
\indfork
\quad\Longrightarrow\quad
\indKim
\quad\Longrightarrow\quad
\indld.
\]
We next show that none of the three reverse implications
holds when $\D$ is proper.
The first non-implication uses the following standard consequence
of the classical uniqueness-of-coheirs characterization of
stability; see
\cite[Theorem~3.10]{BenYaacovStable}.

\begin{fact}
\label{lem:stationarity-coheirs}
Let $T_0$ be a complete theory with monster model $\U$,
and let $\indgen$ be a relation on triples
$M\preccurlyeq A,B\preccurlyeq\U$ of small elementary submodels.
Assume stationarity over models in the sense of
Lemma~\ref{lem:papa-criterion}.
Suppose also that
\[
\operatorname{tp}(A/B)
\text{ is finitely satisfiable in }M
\quad\Longrightarrow\quad
A\indgen_M B.
\]
Then $T_0$ is stable.
\end{fact}

\begin{proof}
Let $M\preccurlyeq B$ and let $q_0,q_1\in S(B)$
be finitely satisfiable in $M$ with the same restriction to $M$.
By extension of coheirs in additional variables, extend them
to types over $B$, finitely satisfiable in $M$, whose common
restriction to $M$ is the type of an
elementary model containing $M$.
The hypothesis and stationarity imply that these extended
types coincide, so $q_0=q_1$.
Since coheirs extend over larger parameter sets, this also
gives uniqueness over arbitrary parameter sets containing $M$.
The characterization of stability cited above completes the proof.
\end{proof}

\begin{theorem}
\label{thm:residual-kim-comparison}
If $\D$ is proper, the following
reverse implications fail on triples
$M\preccurlyeq A,B\preccurlyeq\U$:
\[
\ind
\quad\centernot{\Longleftarrow}\quad
\indfork
\quad\centernot{\Longleftarrow}\quad
\indKim
\quad\centernot{\Longleftarrow}\quad
\indld.
\]
\end{theorem}

\begin{proof}
The theory $T$ is NSOP$_1$ by the Frobenius-field result
recalled in Section~\ref{sec:realization}.
Since $\D$ is proper, it is not simple by
Corollary~\ref{cor:unbounded}.
In particular, $T$ is unstable.

\emph{The first converse.}
The relation $\ind$ has stationarity over models by
Proposition~\ref{prop:stat}.
Since $T$ is unstable, the contrapositive of
Fact~\ref{lem:stationarity-coheirs} gives
$M\preccurlyeq A,B\preccurlyeq\U$ such that
$\operatorname{tp}(A/B)$ is finitely satisfiable
in $M$, but $A\nind_M B$.
A type finitely satisfiable in $M$ extends to a global
coheir, hence to a global $M$-invariant type.
It therefore does not fork over $M$
\cite[Fact~7.2(3)]{KRKim}.
Thus $A\indfork_M B$ but $A\nind_M B$,
with $\operatorname{tp}(A/B)$ finitely satisfiable in $M$.

\emph{The middle converse.}
Since $T$ is NSOP$_1$ and not simple, Kaplan--Ramsey's
criterion
\cite[Proposition~8.4]{KRKim}
gives a model $M$ and tuples $a,b$ such that
\[
a\indKim_M b,
\qquad
a\nindfork_M b.
\]
We include the passage from tuples to elementary models.
Choose $B_0\preccurlyeq\U$ containing $Mb$.
By extension for Kim-independence, choose
$a'\equiv_{Mb}a$ with $a'\indKim_M B_0$.
An automorphism over $Mb$ taking $a'$ to $a$ sends $B_0$
to a model $B$ containing $Mb$, with $a\indKim_M B$.

Choose $A_0\preccurlyeq\U$ containing $Ma$.
By symmetry and extension, an enumerated copy
$B'\equiv_{Ma}B$ can be chosen with $B'\indKim_M A_0$.
An automorphism over $Ma$ taking this enumeration of $B'$
to the corresponding enumeration of $B$ sends $A_0$
to a model $A$ containing $Ma$.
Thus $B\indKim_M A$, and symmetry gives $A\indKim_M B$.
If $A\indfork_M B$, monotonicity of ordinary nonforking
would give $a\indfork_M b$, a contradiction.

\emph{The last converse.}
Let us choose a nontrivial finite simple group $H\notin\D$.
Since $H$ is simple, $H\in\D^{\perp}$. We fix a small model $M\models T$.
By Theorem~\ref{thm:infinite}, $M$ has infinitely many
$H$-Galois extensions (i.e. Galois extensions with Galois group isomorphic to $H$) in $M^{\operatorname{alg}}$.

\par\medskip\noindent
\textbf{Claim 1.}
\emph{Some elementary extension $A$ of $M$ has an
$H$-Galois extension $L/A$ such that $L/M$ is regular.}

\begin{proof}[Proof of Claim 1]
Choose pairwise distinct $H$-Galois extensions $K_i/M$,
$i\in\mathbb N$, inside $M^{\mathrm{alg}}$.
Since $H$ is simple, $K_i\cap K_j=M$ for $i\ne j$.
Let $\mathcal V$ be a nonprincipal ultrafilter on $\mathbb N$,
and put
\[
A=M^{\mathbb N}/\mathcal V,
\qquad
L=\Bigl(\prod_{i\in\mathbb N}K_i\Bigr)/\mathcal V.
\]
Then $M\preccurlyeq A$ via the diagonal embedding.
After choosing isomorphisms $H\cong\Gal(K_i/M)$, the induced coordinatewise action of $H$ on $L$ is faithful,
and $L^H=A$, since $H$ is finite.
Hence we obtain 
$\Gal(L/A)\cong H$.

An irreducible polynomial $f\in M[X]$ of degree greater than one
has a root in at most one $K_i$: otherwise normality would put
its nontrivial splitting field inside two of the $K_i$,
contradicting their intersection being $M$.
By {\L}o\'s's theorem, $f$ has no root in $L$.
Thus $M$ is relatively algebraically closed in $L$,
and characteristic zero gives that $L/M$ is regular.
\end{proof}
We take two $M$-isomorphic copies $(A_1,L_1)$ and $(A_2,L_2)$
of $(A,L)$ from Claim 1 and define
\[
P:=\Frac(L_1\otimes_M L_2),\qquad E:=A_1A_2.
\]
Regularity of $L_i/M$ ensures that the tensor product
is a domain and that $P/L_i$ is regular.
We identify the fields with their canonical images in $P$;
in particular, $A_1\indld_M A_2$.
We will embed $P$ in $\U^{\operatorname{alg}}$ so that
\[
A_1A_2\neq
\acl_\U(A_1A_2)\cap A_1^{\mathrm{alg}}A_2^{\mathrm{alg}}.
\]
Then Fact~\ref{fact:kim-linear-disjointness} witnesses the
failure of Kim-independence.

The two coordinate actions give $|H|^2$ distinct
$E$-automorphisms of $P$, while
\[
[P:E]\leq [L_1:A_1][L_2:A_2]=|H|^2.
\]
Hence $P/E$ is Galois and
\[
\Gal(P/E)\cong H\times H.
\]

Let us define
\[
\Delta_H:=\{(h,h):h\in H\},
\qquad
C:=P^{\Delta_H}.
\]
Then
\[
[C:E]=|H|>1.
\]
\textbf{Claim~2.}
\emph{The extension $C/A_i$ is regular for $i=1,2$.}

\begin{proof}[Proof of Claim~2]
Fix $i\in\{1,2\}$.
Since $P/L_i$ is regular and $L_i/A_i$ is finite,
the relative algebraic closure of $A_i$ in $P$ is $L_i$.
Since $\Delta_H$ projects surjectively onto each factor,
its restriction to $L_i$ is the full group $\Gal(L_i/A_i)$.
Thus
\[
C\cap L_i=L_i^{\Gal(L_i/A_i)}=A_i.
\]
It follows that $A_i$ is relatively algebraically closed in $C$.
Since the characteristic is zero, $C/A_i$ is regular.
\end{proof}
As $A_i\in\K_\D$, we have $A_i^\D=A_i$.
Lemma~\ref{lem:residual-regular-intersection} and Claim 2 yield
\[
C^\D\cap A_i^{\mathrm{alg}}=A_i,
\]
so the extension $C^\D/A_i$ is regular.

By Fact~\ref{fact:galois-residual}\ref{item:galois-absolute-group}
and Fact~\ref{fact:residual-calculus}\ref{item:residual-idempotent},
\[
\RD(\Gal(C^\D))=\Gal(C^\D).
\]
Lemma~\ref{lem:realize} therefore supplies a regular
extension $F'/C^\D$ with $F'\in\K_\D$.
Then, $F'/A_i$ is regular and 
Proposition~\ref{prop:criteria}(b) gives
$A_i\preccurlyeq F'$.
Thus $F'\models T$, and we may embed $F'$ into $\U$
over $M$. We extend this embedding to an embedding of $(F')^{\operatorname{alg}}$ into $\U^{\operatorname{alg}}$. In particular, $P$ is embedded into $\U^{\operatorname{alg}}$ now as required above. Under this embedding, the images of $L_i$ lie in
$A_i^{\mathrm{alg}}$ and their compositum is the image of $P$.
Since $C/E$ is finite and $C\subseteq\U$, we obtain
\[
A_1A_2=E
\subsetneq C
\subseteq
\acl_\U(A_1A_2)\cap A_1^{\mathrm{alg}}A_2^{\mathrm{alg}},
\]
which gives $A_1\nindKim_M A_2$ as indicated above.
\end{proof}

Proposition~\ref{prop:residual-invariant-extension} and
Theorem~\ref{thm:residual-kim-comparison} establish
part~\ref{item:main-kim} of Theorem~\ref{thm:main},
completing its proof.
\par\bigskip

We finish this section with several remarks on the ternary
relations considered above, their descriptions in the
$\omega$-free case, and the failure of local character for $\ind$.
\begin{remark}
Suppose that $\D$ is proper. 
The failure of the implication from $\indKim$ to $\ind$
also follows from Fact~\ref{lem:stationarity-coheirs}:
if the two relations coincided on triples of models,
every coheir would satisfy $\ind$, and stationarity would
force $T$ to be stable.
\end{remark}

\begin{remark}
Kaplan--Ramsey identify Kim-independence over models in
Frobenius fields with Chatzidakis's weak independence
\cite[Definition~9.31 and Theorem~9.33]{KRKim};
see also \cite[Definition~2.3]{Chatz}.
Their description combines independence in a separably
closed overfield with nonforking in the stable theory
of the associated Galois systems.
In characteristic zero, the field-theoretic condition
is algebraic independence in an algebraically closed
overfield.

Our relation imposes the closure condition
$\acl_\U(AB)=(AB)^\D$.
Theorem~\ref{thm:residual-kim-comparison} shows that it
gives a strict refinement of Kim-independence over
models whenever $\D$ is proper.
\end{remark}

\begin{remark}
\label{rem:omega-free-independence-descriptions}
Suppose that $\D=\{1\}$, so that $F$ is an $\omega$-free
PAC field of characteristic zero and $E^\D=E$ for every
field $E$. Let $M\preccurlyeq A,B\preccurlyeq\U$ be small
elementary submodels. All three relations below require
$A\indld_M B$. Under this common assumption, their additional
conditions are
\[
\begin{aligned}
A\ind_M B
&\quad\Longleftrightarrow\quad
(AB)_\U^{\mathrm{alg}}=AB,
\\[4pt]
A\indfork_M B
&\quad\Longleftrightarrow\quad
\begin{gathered}[t]
(AB)_\U^{\mathrm{alg}}
\cap (AB_0)^{\mathrm{alg}}B^{\mathrm{alg}}
=
(AB_0)_\U^{\mathrm{alg}}B
\\
\text{for every }B_0\subseteq B
\text{ with }(B_0)_\U^{\mathrm{alg}}=B_0,
\end{gathered}
\\[4pt]
A\indKim_M B
&\quad\Longleftrightarrow\quad
(AB)_\U^{\mathrm{alg}}
\cap A^{\mathrm{alg}}B^{\mathrm{alg}}=AB.
\end{aligned}
\]
The first equivalence is our definition. The second is
Chatzidakis's description of ordinary nonforking in
$\omega$-free PAC fields \cite{ChatzOmegaFree}; see also
\cite[Example~5.5]{DobrowolskiHoffmannRank} for this explicit
formulation. The third follows from Kaplan--Ramsey's
identification of Kim-independence over models with weak
independence \cite[Theorem~9.33]{KRKim}, together with the description of
weak independence in \cite{ChatzOmegaFree}.

Thus $\ind$ requires the whole compositum $AB$ to be
relatively algebraically closed in $\U$, whereas $\indKim$
only requires that no element of
$(AB)_\U^{\mathrm{alg}}\setminus AB$ lie in
$A^{\mathrm{alg}}B^{\mathrm{alg}}$.
Ordinary nonforking imposes the displayed family of
intersection identities; taking $B_0=M$ recovers the
identity for Kim-independence.
By Theorem~\ref{thm:residual-kim-comparison}, both implications
in
\[
\ind\ \Longrightarrow\ \indfork\ \Longrightarrow\ \indKim
\]
are strict, even on triples of elementary submodels.
\end{remark}

\begin{remark}[Failure of local character]
\label{rem:residual-local-character}
The relation $\ind$ admits a natural extension beyond
elementary submodels. For relatively algebraically closed
subfields $C\subseteq A\cap B$ of $\U$, put
\[
A\ind_C B
\quad\Longleftrightarrow\quad
A\indld_C B
\ \text{ and }\
\acl_\U(AB)=(AB)^\D.
\]
For arbitrary small sets $X,Y,Z$, define $X\ind_Z Y$ by
applying this condition to
\[
C=\acl_\U(Z),\qquad
A=\acl_\U(X\cup Z),\qquad
B=\acl_\U(Y\cup Z).
\]
This extended relation satisfies invariance, symmetry, existence,
extension, monotonicity, base monotonicity, transitivity,
finite character, and stationarity over algebraically closed
bases. We briefly indicate the verifications.
The proofs of Propositions~\ref{prop:exist} and~\ref{prop:stat}
apply to relatively algebraically closed subfields:
such a field $A$ satisfies $A^\D=A$, and
Lemma~\ref{lem:realize} and Proposition~\ref{prop:criteria}(c)
give the required realizations and partial elementary maps.
Monotonicity follows from
Lemma~\ref{lem:residual-regular-intersection};
base monotonicity and transitivity follow by regular base
change and
Fact~\ref{fact:residual-calculus}\ref{item:residual-normal}.
Finite character follows by descending finite Galois
extensions along directed unions of the relevant regular
composita. Extension then follows from existence,
monotonicity and stationarity.

If $\D$ is proper, however, local character fails.
Otherwise, the properties established above, including
stationarity over models, would imply that $T$ is stable
by the axiomatic characterization of stability recalled in
\cite[Section~1, discussion preceding Definition~1.3]{MartinPizarroInd}.
This contradicts Corollary~\ref{cor:unbounded},
since every stable theory is simple.
\end{remark}

\section{Existential closedness for fields with group actions}
\label{sec:regular-ec}
\label{subsec:arbitrary-group-closure}

We study existentially closed and pseudo-existentially closed
fields with group actions.
We begin with the general closure notions and the ultraproduct
criterion for e.c.\ models.
Subsection~\ref{subsec:root-predicate-automorphisms} applies
Corollary~\ref{cor:noMC} to rule out a model companion for
the theory $T_{S,\sigma}$ introduced in
Section~\ref{sec:realization}.
We then compare p.e.c.\ difference fields with p.e.c.\ fields
carrying finite group actions in
Subsections~\ref{subsec:pec-difference-fields}
and~\ref{sec:finite-group-actions}.

Let $G$ be an arbitrary group, and use the language
\[
L_G=\{0,1,+,-,\cdot\}\cup\{\sigma_g:g\in G\}.
\]
A \emph{$G$-field} is a field with an action of $G$ by
automorphisms; faithfulness is not part of this definition.
We initially allow every characteristic, including imperfect
fields. A field extension is \emph{regular} if it is separable
and the base field is relatively algebraically closed in it.
The following terminology refers to the ordinary language $L_G$;
see \cite[Definition~2.1 and Remark~2.3]{BKTorsion}.
In particular, \emph{p.e.c.} always refers to this ordinary-language
notion.

\begin{definition}\label{def:group-closure-notions}
A $G$-field $K$ is
\begin{enumerate}[label=\textup{(\alph*)}]
\item
\emph{existentially closed}, or \emph{e.c.}, if it is
existentially closed in every $G$-field extension;
\item
\emph{pseudo-existentially closed}, or \emph{p.e.c.}, if it is
existentially closed in every $G$-field extension whose underlying
field extension is regular;
\item
\emph{algebraically $G$-closed}, or simply \emph{$G$-closed},
if it has no proper algebraic $G$-field extension.
\end{enumerate}
In the last condition, algebraic extensions may be infinite.
\end{definition}

For the trivial group, p.e.c. fields are exactly PAC fields
\cite[Proposition~11.3.5]{FJ} and e.c. fields are exactly algebraically closed fields.
For $G=\mathbb Z$, a $G$-field is a difference field with a
distinguished automorphism.

The following observation is recorded in
\cite[Remark~2.3(0)--(1)]{BKTorsion}. Neither finiteness nor
finite generation of $G$ is needed.

\begin{fact}
\label{prop:ec-pec-gclosed}
For every group $G$ and every $G$-field $K$,
\[
K\text{ is e.c.}
\quad\Longleftrightarrow\quad
K\text{ is p.e.c. and algebraically $G$-closed}.
\]
\end{fact}

For a first-order theory $T_0$, let $\operatorname{EC}(T_0)$
denote the class of all models of $T_0$ that are existentially
closed in every extension which is also a model of $T_0$.
Here \emph{elementary} means axiomatizable by a possibly
incomplete first-order theory.
We recall the following standard facts; see also
\cite[Introduction and discussion following Fact~5.3]{GSTEcActions}
for a recent discussion of the ultraproduct criterion.

\begin{fact}
\label{fact:ec-ultraproducts}
Let $T_0$ be a first-order theory.
\begin{enumerate}[label=\textup{(\alph*)}]
\item
The class $\operatorname{EC}(T_0)$ is closed under
elementary substructures.
\item
The class $\operatorname{EC}(T_0)$ is elementary if and only
if it is closed under ultraproducts.
\item
If $T_0$ is consistent and inductive, the equivalent
conditions in \textup{(b)} are also equivalent to the
existence of a model companion of $T_0$.
\end{enumerate}
\end{fact}

\begin{proof}
Part~\textup{(a)} follows from the elementary-substructure
observation in \cite[Section~3.2]{TentZiegler}.

For~\textup{(b)}, part~\textup{(a)} gives closure under
ultraroots: if an ultrapower of $M$ belongs to
$\operatorname{EC}(T_0)$, then so does $M$, since the
diagonal embedding is elementary.
As $\operatorname{EC}(T_0)$ is also closed under
isomorphisms, the assertion follows from the
ultraproduct--ultraroot criterion
\cite[Corollary~16.2.6]{GoldbringUltrafilters}.

Part~\textup{(c)} is
\cite[Proposition~3.5.15]{ChangKeisler}.
\end{proof}

Parts~\textup{(a)} and~\textup{(b)} require no completeness,
inductiveness, or amalgamation hypothesis on $T_0$.
Inductiveness is used only in~\textup{(c)}.

For an inductive theory, nonexistence of a model companion
therefore already implies failure of closure of its e.c.\
class under ultraproducts.
This gives an ultraproduct of possibly different e.c.\
models which is not e.c.
Producing such a witness with all factors equal is a
stronger assertion and requires an additional argument.
The orbit construction below provides a witness of this
more precise kind.

The criterion concerns existential closure in a fixed
first-order language with respect to all embeddings between
models of a fixed theory.
It cannot be applied to p.e.c.\ difference fields merely
by restricting the permitted extensions to regular ones.
Although root predicates make embeddings correspond to
regular extensions, they also change which formulas are
existential; see Remark~\ref{rem:root-ec-distinction}.

We next introduce the languages used for difference fields,
that is, for $G=\mathbb Z$, with a generator acting as $\sigma$.
The ordinary difference language is
\[
L_\sigma=\{0,1,+,-,\cdot,\sigma\},
\]
where $\sigma$ is required to be an automorphism.
Naming all positive and negative iterates gives $L_{\mathbb Z}$.
Inverse iterates can be eliminated from existential formulas by
introducing variables and equations of the form $\sigma(y)=x$,
so the p.e.c.\ notions in the two languages agree.

\subsection{Existentially closed models in the root-predicate language}
\label{subsec:root-predicate-automorphisms}

We first study e.c.\ models of the theory $T_{S,\sigma}$
from Section~\ref{sec:realization}, and deduce the nonexistence
of a model companion from Corollary~\ref{cor:noMC}.
We begin by distinguishing existential closedness in the
expanded language from the ordinary-language p.e.c.\ condition.

\begin{remark}\label{rem:root-ec-distinction}
The allowable extensions are the same as in the definition of
p.e.c., but the expanded language has additional existential
formulas. For example,
\[
\exists a\,
\bigl(\sigma(a)=a\ \land\ \neg S_2(0,-a)\bigr)
\]
is existential in $\mathcal L_{S,\sigma}$, whereas its defining
translation into $L_\sigma$ is
\[
\exists a\,
\bigl(\sigma(a)=a\ \land\ \forall x\,(x^2\ne a)\bigr).
\]

For any difference field $(K,\sigma)$, extend $\sigma$ to
$K(t)$ by fixing a transcendental element $t$.
The extension $K(t)/K$ is regular, and the displayed
expanded-language formula holds with witness $a=t$.
It follows that every e.c.\ model of $T_{S,\sigma}$ contains
an element fixed by $\sigma$ that is a nonsquare in its
underlying field.

On the other hand, a model $(K,\sigma)$ of
$\mathrm{ACFA}_0$ is e.c.\ in the ordinary difference language,
and hence p.e.c., but its underlying field is algebraically
closed. The displayed formula fails in $K$, so
$(K,\sigma)_S$ is not e.c.\ as a model of $T_{S,\sigma}$.

Consequently, the $L_\sigma$-reducts of the e.c.\ models of
$T_{S,\sigma}$ form a proper subclass of the p.e.c.\
characteristic-zero difference fields.
Inclusion follows from
Proposition~\ref{prop:root-regular-embeddings}, and the
$\mathrm{ACFA}_0$ example shows that it is strict.
Indeed, none of those reducts is a model of $\mathrm{ACFA}_0$.
\end{remark}

Let $T_\omega$ be the theory of characteristic-zero
$\omega$-free PAC fields in the canonical root-predicate language.
It is the model companion of the theory of characteristic-zero
fields in that language \cite[Theorem~27.2.3]{FJ}.

\begin{proposition}
\label{prop:root-ec-underlying}
If $(K,\sigma)_S$ is an e.c.\ model of $T_{S,\sigma}$,
then $K$ is an $\omega$-free PAC field.
\end{proposition}

\begin{proof}
Let $E/K$ be a regular field extension.
Take independent copies $E_i$ of $E$, indexed by $i\in\mathbb Z$,
with isomorphisms $\iota_i:E\longrightarrow E_i$ satisfying
$\iota_i|_K=\sigma^i$, and take $E_0=E$.
Let $L$ be the fraction field of their tensor product over $K$.
Regularity of $E/K$ implies that this tensor product is a domain
and that $L/E_i$ is regular for every $i$.
The shift
\[
\iota_i(e)\longmapsto\iota_{i+1}(e)
\]
defines an automorphism $\tau$ of $L$ extending $\sigma$.
Since $L/E$ is regular, every existential $\mathcal L_S$-formula
over $K$ realized in $E_S$ is also realized in $L_S$.
Existential closedness of $(K,\sigma)_S$ reflects it back to $K_S$.
Thus $K_S$ is existentially closed among characteristic-zero
fields with their canonical root-predicate expansions.
By \cite[Theorem~27.2.3]{FJ}, these existentially closed
structures are precisely the canonical expansions of
$\omega$-free PAC fields. Hence $K$ is an $\omega$-free
PAC field, as required.
\end{proof}

\begin{theorem}\label{thm:root-automorphism-no-companion}
The theory $T_{S,\sigma}$ has no model companion.
\end{theorem}

\begin{proof}
Suppose that $T^*$ is a model companion of $T_{S,\sigma}$,
and choose a completion $T$ of $T_\omega$.
Since $T_{S,\sigma}$ is inductive by
Fact~\ref{lem:root-ec-hull}, the models of $T^*$ are precisely
its e.c.\ models. Their $\mathcal L_S$-reducts therefore
satisfy $T_\omega$ by Proposition~\ref{prop:root-ec-underlying}.

Every model $(K,\sigma)_S$ of $T_\sigma$ embeds in a model
$(M,\tau)_S$ of $T^*$.
Model completeness of $T_\omega$ gives $K_S\preccurlyeq M_S$,
so $M_S\models T$.
Consequently, $T^*\cup T$ is a model companion of $T_\sigma$:
it is model complete, its models satisfy $T_\sigma$, and
every model of $T_\sigma$ embeds in one of its models.
For $\D=\{1\}$, the class $\K_\D$ consists precisely of the
characteristic-zero $\omega$-free PAC fields.
Corollary~\ref{cor:noMC} therefore rules out such a model
companion of $T_\sigma$, a contradiction.
\end{proof}

\begin{remark}\label{rem:root-ec-ultraproducts}
Since $T_{S,\sigma}$ is inductive by
Fact~\ref{lem:root-ec-hull},
Theorem~\ref{thm:root-automorphism-no-companion} and
Fact~\ref{fact:ec-ultraproducts}\textup{(c)}
imply that its e.c.\ models are not closed under ultraproducts.

In fact, closure under ultrapowers also fails.
We omit the proof, which follows the orbit construction in
the proof of Theorem~\ref{thm:regular-ec-nonelementary} below,
with the initial regular extension of $(\mathbb Q(z),\sigma)$,
$\sigma(z)=z+1$, chosen e.c.\ in the expanded language using
Fact~\ref{lem:root-ec-hull}.
Proposition~\ref{prop:root-regular-embeddings} ensures that
the regular extension used in that construction is admissible
in the expanded language.
\end{remark}

\subsection{P.e.c. difference fields}
\label{subsec:pec-difference-fields}

We now work in the ordinary difference language $L_\sigma$.
Throughout this subsection the characteristic is zero, so
regularity means that the base field is relatively algebraically
closed in the extension.
We use the preservation of regularity under subextensions over
the fixed base, towers, and directed unions.
The p.e.c.\ notion for difference fields is also discussed in
\cite[Remark~4.9 and Question~4.10]{BK}.

Taking $L_\sigma$-reducts of the e.c.\ extensions supplied by
Fact~\ref{lem:root-ec-hull} gives the following consequence.

\begin{fact}
\label{lem:regular-ec-hull}
Every characteristic-zero difference field $(K,\sigma)$ has a
regular difference-field extension $(M,\tau)$ that is p.e.c.
\end{fact}

\begin{proof}
Apply Fact~\ref{lem:root-ec-hull} to the canonical expansion
$(K,\sigma)_S$, obtaining an e.c. $T_{S,\sigma}$-extension
$(M,\tau)_S$ with $M/K$ regular.
By Proposition~\ref{prop:root-regular-embeddings}, every regular
difference-field extension of $(M,\tau)$ canonically expands
to a $T_{S,\sigma}$-extension of $(M,\tau)_S$.
Every existential formula in the ordinary difference language
is also existential in the expanded language.
Consequently, existential closure of $(M,\tau)_S$ implies that
$(M,\tau)$ is p.e.c.
\end{proof}

\begin{lemma}\label{lem:disjoint-orbits}
Let $(M,\sigma)$ be a characteristic-zero difference field, and let
$a,b\in M$ have disjoint infinite orbits under $\sigma$. Then there
is a regular difference-field extension $(L,\sigma)/(M,\sigma)$
containing $t,u,v$ such that
\[
\sigma(t)=t,\qquad u^2=t-a,\qquad v^2=b-t,\qquad uv\ne0.
\]
\end{lemma}

\begin{proof}
Write $\alpha_i=\sigma^i(a)$ and $\beta_i=\sigma^i(b)$ for
$i\in\mathbb Z$. By hypothesis all elements in these two indexed
families are pairwise distinct. Take $t$ transcendental over $M$ and
form
\[
L=M(t)(u_i,v_i:i\in\mathbb Z),\qquad
u_i^2=t-\alpha_i,\quad v_i^2=\beta_i-t.
\]
This is a standard construction, see e.g. the proof of \cite[Lemma 3.1]{BK2026}, and the extension $L/M$ is regular.

We extend $\sigma$ to $L$ by setting
\[
\sigma(t)=t,\qquad \sigma(u_i)=u_{i+1},\qquad
\sigma(v_i)=v_{i+1}.
\]
The elements $t,u_0,v_0$ have the
required properties; the two roots are nonzero because $t$ is
transcendental over $M$.
\end{proof}

\begin{theorem}\label{thm:regular-ec-nonelementary}
The class of p.e.c. difference fields of
characteristic zero is not closed under ultrapowers. In particular,
it is not elementary.
\end{theorem}

\begin{proof}
Let $a$ be transcendental and $\sigma\in \operatorname{Aut}(\mathbb Q(a))$ such that
$\sigma(a)=a+1$. By Fact~\ref{lem:regular-ec-hull}, there is a
p.e.c. regular difference-field extension
$(K,\sigma)/(\mathbb Q(a),\sigma)$. Since $-1$ is not a square in $\mathbb Q(a)$ and
$\mathbb Q(a)$ is relatively algebraically closed in $K$,
it is not a square in $K$ either. We consider the following existential formula
\[
\psi(y,z):=\exists x\,\exists u\,\exists v\,
\bigl(\sigma(x)=x\ \land\ u^2=x-y\ \land\ v^2=z-x
\ \land\ uv\ne0\bigr).
\]
For every positive integer $n$, one has
\[
K\models\neg\psi(a,\sigma^n(a)).
\]
Otherwise, applying $\sigma^n$
to $u^2=x-a$ gives
$(\sigma^n(u))^2=x-\sigma^n(a)=-v^2$.
Since $u$ and $v$ are nonzero, this would make $-1$ a square in $K$.

Let $\mathcal V$ be a nonprincipal ultrafilter on the positive
integers and let $M=K^{\mathbb N_{>0}}/\mathcal V$.
In $M$, write $a$ also for the constant class of $a$, and put
\[
b=[(\sigma^n(a))_{n>0}]_{\mathcal V}.
\]
By \L o\'s's theorem, $M\models\neg\psi(a,b)$.
Nevertheless, $a$ and $b$ have disjoint infinite orbits in $M$.
Indeed, their iterates are $a+i$ and $b+j$ for $i,j\in\mathbb Z$.
Each orbit has distinct elements, and an equality between its two
displayed forms would require
\[
\{n>0:a+i=a+n+j\}\in\mathcal V.
\]
This set contains at most one integer, which is impossible.

Lemma~\ref{lem:disjoint-orbits} now supplies a regular
difference-field extension of $M$ satisfying $\psi(a,b)$.
Thus $M$ is not p.e.c., although $K$ is.
\end{proof}

\begin{remark}\label{rem:pec-versus-no-companion}
The non-elementarity assertion of
Theorem~\ref{thm:regular-ec-nonelementary}
does not follow formally from either
Corollary~\ref{cor:noMC} or
Theorem~\ref{thm:root-automorphism-no-companion}.
Those results concern existential closure in the
root-predicate language.
As Remark~\ref{rem:root-ec-distinction} shows, the ordinary
reducts of its e.c.\ models form a proper subclass of the
p.e.c.\ difference fields.
Non-elementarity of that subclass does not imply
non-elementarity of the entire p.e.c.\ class.
\end{remark}

\subsection{Finite group actions}
\label{sec:finite-group-actions}
\label{subsec:finite-group-characterization}

We now assume that $G$ is finite and again allow arbitrary
characteristic, including imperfect fields.
This subsection concerns p.e.c. in the ordinary language $L_G$.
The elementarity of this class is recorded in
\cite[Remark~6.25]{Gogolok}.
We give an explicit characterization in terms of faithfulness
and the PAC property of the fixed field.

\begin{theorem}\label{thm:finite-group-pec}
A $G$-field $K$ is p.e.c. if and only if the action of $G$
is faithful and the fixed field $C=K^G$ is PAC.
\end{theorem}

\begin{proof}
Suppose first that $K$ is p.e.c. Adjoin algebraically independent
variables $(t_h)_{h\in G}$ and extend the action to
$K(t_h:h\in G)$ by $\sigma_g(t_h)=t_{gh}$.
This is a regular $G$-field extension and, for every $g\ne1$,
it satisfies $\exists x\,(\sigma_g(x)\ne x)$.
The same sentence holds in $K$, so the action is faithful.
Artin's fixed-field theorem gives
$\Gal(K/C)\cong G$ and $[K:C]=|G|$.

Let $E/C$ be a regular field extension.
The tensor product $K\otimes_C E$ is a domain and is finite
dimensional over $E$, so it is a field $L$.
The extension $L/K$ is regular by base change.
Extend the action of $G$ by fixing $E$ pointwise.
An existential field formula over $C$ realized in $E$ is
therefore realized in $L$ by coordinates fixed by every
element of $G$. Adjoin these finitely many fixed-point
equations to the formula. The p.e.c. property gives such
a witness in $K$, and its coordinates lie in $C$.
Thus $C$ is existentially closed in every regular extension,
and hence is PAC.

Conversely, suppose that the action is faithful and $C$ is PAC.
Let $L/K$ be a regular $G$-field extension and put $D=L^G$.
The action on $L$ is faithful, and Artin's theorem gives
\[
\Gal(K/C)\cong G\cong\Gal(L/D),\qquad
[K:C]=[L:D]=|G|.
\]
The subgroup of $\Gal(L/D)$ fixing $KD$ is trivial, because
the action on $K$ is faithful. Hence $L=KD$.
It follows from the equality of degrees that a $C$-basis
of $K$ remains a $D$-basis of $L$.

The extension $D/C$ is regular. Indeed, $K/C$ is finite
separable and $L/K$ is separable, so $L/C$, and therefore
$D/C$, is separable. If $d\in D$ is algebraic over $C$,
then it is algebraic over $K$. Regularity of $L/K$ gives
$d\in K$, and therefore $d\in K\cap D=C$.

Choose a $C$-basis $e_1,\ldots,e_m$ of $K$.
Write each element of $L$ uniquely as $\sum_i e_i y_i$,
with $y_i\in D$.
Multiplication and the maps $\sigma_g$ have structure
constants in $C$ in this basis; the same formulas describe
their actions on $L$, since $G$ fixes $D$ pointwise.
Parameters from $K$ also have coordinates in $C$.
Consequently, an existential $L_G$-formula over $K$
translates into an existential field formula over $C$:
equality becomes equality of all coordinates, and inequality
becomes a disjunction of coordinate inequalities.
A witness in $L$ gives a coordinate witness in $D$.
Since $C$ is PAC and $D/C$ is regular, there is a coordinate
witness in $C$, which gives a witness to the original formula
in $K$. Thus $K$ is p.e.c.
\end{proof}

\begin{corollary}\label{cor:finite-group-elementary}
For each fixed finite group $G$, each of the following classes
of $G$-fields is elementary:
\begin{enumerate}[label=\textup{(\alph*)}]
\item the p.e.c. $G$-fields;
\item the algebraically $G$-closed $G$-fields;
\item the e.c. $G$-fields.
\end{enumerate}
These assertions hold without a restriction on characteristic.
\end{corollary}

\begin{proof}
For the first assertion, add to the field and action axioms
the faithfulness sentences
\[
\exists x\,(\sigma_g(x)\ne x)
\qquad(g\in G\setminus\{1\}).
\]
The fixed field is defined by the single formula
\[
\operatorname{Fix}(x):=
\bigwedge_{g\in G}\sigma_g(x)=x.
\]
Relativize the PAC axiom scheme to this formula.
Theorem~\ref{thm:finite-group-pec} identifies the models of
the resulting theory with the p.e.c. $G$-fields.

For the second assertion, any proper algebraic $G$-field
extension $L/K$ contains a proper finite $G$-field extension:
for $a\in L\setminus K$ algebraic over $K$, take
$K(\sigma_g(a):g\in G)$.
For each $n\geqslant2$, the existence of an extension of
degree $n$ carrying the prescribed $G$-action is first-order
expressible over $K$. One codes an $n$-dimensional field
extension by multiplication constants on $K^n$, taking
$e_1=(1,0,\ldots,0)$ as its unit and $Ke_1$ as its scalar
copy of $K$. Each $\sigma_g$ is coded by a semilinear
matrix, using the prescribed $\sigma_g$ on scalar
coefficients and fixing $e_1$. The field axioms,
compatibility with multiplication, and the finitely many
group-action identities are first-order conditions.
Negating these sentences for all $n\geqslant2$
axiomatizes algebraic $G$-closedness.

The third assertion follows from the first two and
Fact~\ref{prop:ec-pec-gclosed}.
\end{proof}
The theory of e.c. $G$-fields is the model companion
$G\text{-}\mathrm{TCF}$ studied in \cite{HKFinite}.

\begin{remark}[A p.e.c. field that is not e.c.]
\label{rem:finite-pec-not-ec}
Let $C$ be a pseudofinite field of characteristic zero.
Then $C$ is PAC and $\Gal(C)\cong\widehat{\mathbb Z}$.
Choose a cyclic extension $L/C$ of degree $6$, let $K/C$
be its quadratic subextension, and let $\tau$ generate
$\Gal(L/C)$.
Equip $K$ with the action of $G=\mathbb Z/2\mathbb Z$
whose nonidentity element acts as $\tau^3|_K$.
This action is faithful and $K^G=C$, so
Theorem~\ref{thm:finite-group-pec} makes $K$ p.e.c.
The automorphism $\tau^3$ on $L$ has order $2$ and extends
the given action. Hence $L/K$ is a proper algebraic
$G$-field extension, and $K$ is not e.c.
\end{remark}

\subsection{General conclusions}
\label{subsec:group-action-conclusions}
\label{subsec:cyclic-closure-contrast}

The decomposition in Fact~\ref{prop:ec-pec-gclosed}
holds for every group. Its consequences for axiomatizability
differ between finite groups and $\mathbb Z$.
The following observation about $\mathbb Z$-closedness
holds in every characteristic.

\begin{proposition}\label{prop:z-closed-acf}
A difference field $(K,\sigma)$ is algebraically
$\mathbb Z$-closed if and only if $K$ is algebraically closed.
Consequently, in every characteristic,
\[
(K,\sigma)\models\mathrm{ACFA}
\quad\Longleftrightarrow\quad
\bigl((K,\sigma)\text{ is p.e.c.}\bigr)
\ \text{and}\ K\models\mathrm{ACF}.
\]
\end{proposition}

\begin{proof}
Every automorphism $\sigma$ of $K$ extends to an
automorphism $\widetilde\sigma$ of an algebraic closure
$K^{\mathrm{alg}}$. Thus
$(K^{\mathrm{alg}},\widetilde\sigma)/(K,\sigma)$ is an
algebraic $\mathbb Z$-field extension, and algebraic
$\mathbb Z$-closedness forces $K=K^{\mathrm{alg}}$.
The converse is immediate.
The final equivalence follows from
Fact~\ref{prop:ec-pec-gclosed} and the classical
axiomatization of e.c. difference fields by $\mathrm{ACFA}$;
see \cite{CHDifference}.
\end{proof}

For clarity, write $\mathcal P_{G,0}$, $\mathcal A_{G,0}$,
and $\mathcal E_{G,0}$ for the classes of characteristic-zero
$G$-fields that are respectively p.e.c., algebraically
$G$-closed, and e.c. Fact~\ref{prop:ec-pec-gclosed}
gives, for every group $G$,
\[
\mathcal E_{G,0}=\mathcal P_{G,0}\cap\mathcal A_{G,0}.
\]
The elementary-class behavior differs as follows:
\begin{center}
\begin{tabular}{lcc}
\hline
Class & $G$ finite & $G=\mathbb Z$\\
\hline
$\mathcal P_{G,0}$ & elementary & not elementary\\
$\mathcal A_{G,0}$ & elementary & elementary\\
$\mathcal E_{G,0}$ & elementary & elementary\\
\hline
\end{tabular}
\end{center}
The finite-group column follows from
Corollary~\ref{cor:finite-group-elementary}.
For $G=\mathbb Z$, the last two entries follow from
Proposition~\ref{prop:z-closed-acf} and the existence of
$\mathrm{ACFA}_0$. The first entry is the separate
non-elementarity result of
Theorem~\ref{thm:regular-ec-nonelementary}.
In particular,
\[
\mathcal P_{\mathbb Z,0}
\cap
\{(K,\sigma):K\models\mathrm{ACF}_0\}
=\operatorname{Mod}(\mathrm{ACFA}_0).
\]
Thus p.e.c. is axiomatizable on algebraically closed
underlying fields, although it is not axiomatizable among
all characteristic-zero difference fields.
An elementary intersection with an elementary class does
not force the other class to be elementary.

\begin{corollary}\label{cor:model-companion-not-pec}
Existence of a model companion for the theory of $G$-fields,
even together with axiomatizability of algebraic
$G$-closedness, does not imply axiomatizability of the
p.e.c. $G$-fields. This failure already occurs for
$G=\mathbb Z$ in characteristic zero.
\end{corollary}

\begin{proof}
For $G=\mathbb Z$, the model companion is
$\mathrm{ACFA}_0$, algebraic $G$-closedness is
axiomatized by $\mathrm{ACF}_0$ together with the
automorphism axioms, and
Theorem~\ref{thm:regular-ec-nonelementary} gives
non-elementarity of the p.e.c. class.
\end{proof}

\noindent\textbf{Use of AI.}
The authors used ChatGPT (GPT-6.0 with reasoning
effort set to Ultra) to assist in developing proofs and drafting
portions of the text. The authors refined and streamlined the
arguments and assume full responsibility for the final content.

\end{document}